\documentclass{CustomArticle}

\usepackage{algorithm}
\usepackage[noend]{algpseudocode}
\usepackage[
  style = trad-abbrv,
  citestyle = numeric-comp,
  sorting = nyt
]{biblatex}
\usepackage{booktabs}
\usepackage{caption}
\usepackage{graphicx}
\usepackage{import}
\usepackage{multirow}
\usepackage{hyperref}
\usepackage[capitalize]{cleveref}

\usepackage{CommonMath}
\usepackage{ProjectHelperMacros}
\usepackage{ProjectStyleTweaks}
\usepackage{ProjectFixes}

\PaperGeneralInfo{
  title = {
    Universality of AdaGrad Stepsizes\\
    for Stochastic Optimization:\\
    Inexact Oracle, Acceleration and Variance Reduction
  },
  date = {June 10, 2024}
}
\AddPaperAuthor{
  name = {Anton Rodomanov},
  affiliation = {CISPA Helmholtz Center for Information Security},
  email = {anton.rodomanov@cispa.de}
}
\AddPaperAuthor{
  name = {Xiaowen Jiang},
  affiliation = {CISPA Helmholtz Center for Information Security},
  email = {xiaowen.jiang@cispa.de}
}
\AddPaperAuthor{
  name = {Sebastian Stich},
  affiliation = {CISPA Helmholtz Center for Information Security},
  email = {stich@cispa.de}
}

\PaperAbstract{
  We present adaptive gradient methods (both basic and accelerated) for
  solving convex composite optimization problems in which the main part is
  approximately smooth (a.k.a.\ $(\delta, L)$-smooth) and can be accessed
  only via a (potentially biased) stochastic gradient oracle.
  This setting covers many interesting examples including Hölder smooth
  problems and various inexact computations of the stochastic gradient.
  Our methods use AdaGrad stepsizes and are adaptive in the sense that they
  do not require knowing any problem-dependent constants except an estimate
  of the diameter of the feasible set but nevertheless achieve the best
  possible convergence rates as if they knew the corresponding constants.
  We demonstrate that AdaGrad stepsizes work in a variety of situations by
  proving, in a unified manner, three types of new results.
  First, we establish efficiency guarantees for our methods in the classical
  setting where the oracle's variance is uniformly bounded.
  We then show that, under more refined assumptions on the variance, the same
  methods without any modifications enjoy implicit variance reduction
  properties allowing us to express their complexity estimates in terms of
  the variance only at the minimizer.
  Finally, we show how to incorporate explicit SVRG-type variance reduction
  into our methods and obtain even faster algorithms.
  In all three cases, we present both basic and accelerated algorithms
  achieving state-of-the-art complexity bounds.
  As a direct corollary of our results, we obtain universal stochastic
  gradient methods for Hölder smooth problems which can be used in all
  situations.
}
\PaperKeywords{%
  convex optimization, stochastic optimization, adaptive methods, universal
  algorithms, complexity estimates, acceleration, variance reduction, AdaGrad,
  SVRG, weakly smooth functions, Hölder condition, inexact oracle
}

\begin{document}
  \PrintTitleAndAbstract

  \section{Introduction}

\paragraph{Motivation.}

Gradient methods are the most popular and efficient optimization algorithms
for solving machine learning problems.
To achieve the best convergence speed for these algorithms, their stepsizes
needs to be chosen properly.
While there exist various theoretical recommendations, dictated by the
convergence analysis, on how to select stepsizes based on various
problem-dependent parameters, they are usually impractical because the
corresponding constants may be unknown or their worst-case estimates might be
too pessimistic.
Furthermore, every applied problem usually belongs to multiple problem classes
at the same time, and it is not always evident in advance which of them better
suits the concrete problem instance one works with.
For classical optimization algorithms, this problem is typically resolved by
using a line search.
This is a simple yet powerful mechanism which automatically chooses the best
stepsize by checking at each iteration a certain condition involving the
objective value, its gradient, etc.

However, the line-search approach is unsuitable for problems of stochastic
optimization, where gradients are observed with random noise.
For these problems, it is common instead to apply so-called adaptive methods
which set up their stepsizes by simply accumulating on-the-fly certain
information about observed stochastic gradients.
The first such an algorithm,
AdaGrad~\cite{duchi2011adaptive,mcmahan2010adaptive}, was obtained from
theoretical considerations but quickly inspired several other heuristic
methods like RMSProp~\cite{tieleman2012lecture} and Adam~\cite{KingmaB14} that
are now at the forefront of training machine learning models.

Excellent practical performance of adaptive methods on various applied problems
naturally sparked a lot of theoretical interest in these algorithms.
An important observation was done by \textcite{levy2018online} who showed that
AdaGrad possesses a certain universality property, in the sense that it
works for several problem classes simultaneously.
Specifically, they showed that AdaGrad converges both for nonsmooth problems
with bounded gradient and also for smooth problems with Lipschitz gradient,
without needing to know neither the corresponding Lipschitz constants, nor the
oracle's variance but enjoying the rates which are characteristic for
algorithms which have the knowledge of these constants.
They also presented an accelerated version of AdaGrad with similar properties.
Further improvements and generalization of these ideas were considered
in~\cite{kavis2019unixgrad,joulani2020simpler,ene2021adaptive}.

Nonsmooth and smooth problems are the extremes of the more general Hölder class
of problems.
The fact that AdaGrad methods simultaneously work for these two extreme cases
does not seem to be a coincidence and suggests that these algorithms should
work more generally for any problem with intermediate level of smoothness.
Some further confirmations to this were recently provided
in~\cite{orabona2023normalized} although in a rather restricted setting of
deterministic problems and only for the basic AdaGrad method.
The stochastic case and acceleration were constituting an open problem which
was recently resolved in~\cite{rodomanov2024universal} for a slightly modified
AdaGrad stepsize.

All the previously discussed results were proved only for the classical
stochastic optimization setting where the variance of stochastic gradients is
assumed to be uniformly bounded.
In a recent work, \textcite{attia2023sgd} showed that the basic AdaGrad method
for smooth problems works under the more general assumption when the variance
is bounded by a constant plus a multiple of the squared gradient norm.
On a related note, it was also shown recently that AdaGrad stepsizes can be
used inside gradient methods with SVRG-type variance-reduction.
The first such an algorithm was proposed in~\cite{dubois2022svrg}.
The accelerated SVRG method enjoying optimal worst-case oracle complexity
for smooth finite-sum optimization problems was later presented
in~\cite{liu2022adaptive}.

\paragraph{Contributions.}

In this work, we further extend the results mentioned above by demonstrating
that AdaGrad stepsizes are even more universal than was shown previously in the
literature.
Specifically, we consider the composite optimization problem where the main
part is approximately smooth (a.k.a.\ $(\delta, L)$-smooth) and can be accessed
only via a (potentially biased) stochastic gradient oracle.
This setting is more general than typically considered in the literature on
adaptive methods and covers many interesting examples, including smooth,
nonsmooth and, more generally, Hölder smooth problems, problems in which the
objective function is given itself as another optimization problem whose
solution can be computed only approximately, etc.

Our contributions can be summarized as follows:
\begin{enumerate}
  \item
    We start, in \cref{sec:MainAlgorithmsAndStepsizeUpdateRules}, with
    identifying the key property of the AdaGrad stepsizes which allows us to
    apply these stepsizes, in a unified manner, in a variety of situations we
    consider later.
    We also present our two mains algorithms, $\UniversalSgd$ and
    $\UniversalFastSgd$ which are the classical stochastic gradient method (SGD)
    and its accelerated version, respectively, but equipped with AdaGrad
    stepsizes.

  \item
    We then establish, in \cref{sec:UniformlyBoundedVariance}, efficiency
    guarantees for these methods in the classical setting where the oracle's
    variance is assumed to be uniformly bounded.

  \item
    In \cref{sec:ImplicitVarianceReduction}, we complement these results
    by showing that, under additional assumptions that the variance is
    itself approximately smooth w.r.t.\ the objective function,
    the same $\UniversalSgd$ and $\UniversalFastSgd$ without any modifications
    enjoy implicit variance reduction properties allowing us to express their
    complexity estimates in terms of the variance only at the minimizer.

  \item
    Under the additional assumption that one can periodically compute the full
    (inexact) gradient of the objective function, we show, in
    \cref{sec:ExplicitVarianceReductionWithSvrg}, how to incorporate explicit
    SVRG-type variance reduction into our methods, obtaining new
    $\UniversalSvrg$ and $\UniversalFastSvrg$ algorithms which enjoy even faster
    convergence rates by completely eliminating the variance.
\end{enumerate}

Our results are summarized in \cref{tab:SummaryOfMainResults}
(in the BigO-notation).
In all the situations, we present both basic and accelerated algorithms whose
only essential parameter is an estimate~$D$ of the diameter of the feasible
set;
the methods automatically adapt to all other problem-dependent constants.
In a number of special cases, our algorithms achieve known state-of-the-art
complexity bounds, but not restricted to those special cases.
In \cref{sec:ApplicationToHolderSmoothProblems}, we illustrate the significance
of our results by demonstrating that complexities for our methods on stochastic
optimization problems with Hölder smooth components can be obtained as simple
corollaries from our main results.

\begin{table}
  \crefname{algorithm}{Alg.}{Algs.}
  \crefname{theorem}{Thm.}{Thms.}
  \crefname{assumption}{}{}
  \crefcompressoff
  \renewcommand{\creflastconjunction}{, }
  \centering
  \captionsetup{font=small}
  \caption{%
    Summary of main results for solving problem~\eqref{eq:Problem} with our
    methods.
    ``Convergence rate'' is expressed in terms of the expected function residual
    at iteration~$k$ (or $t$, depending on the method).
    ``SO complexity'' denotes the cumulative stochastic-oracle complexity of
    the method since its start and up to iteration~$k$ (or~$t$), which is
    defined as the number of queries to the stochastic oracle~$\hat{g}$;
    for SVRG methods, we assume that querying the (inexact) full-gradient
    oracle~$\bar{g}$ is $n$ times more expensive than~$\hat{g}$, and define the
    SO complexity as $N_{\hat{g}} + n N_{\bar{g}}$, where $N_{\hat{g}}$ and
    $N_{\bar{g}}$ are the number of queries to~$\hat{g}$ and~$\bar{g}$,
    respectively.
  }
  \label{tab:SummaryOfMainResults}
  \scriptsize
  \begin{tabular}{lcccc}
    \toprule
    Method & Convergence rate & SO complexity & Assumptions & Reference \\
    \midrule
    \multirow{2}{*}{$\UniversalSgd$ (\cref{alg:UniversalSgd})}
           &
    $\frac{L_f D^2}{k} + \frac{\sigma D}{\sqrt{k}} + \delta_f$
           &
    \multirow{2}{*}{$k$}
           &
    \cref{as:ApproximateSmoothness,as:BoundedFeasibleSet,as:UniformlyBoundedVariance}
           &
    \cref{th:UniversalSgd}
    \\
           &
    $
      \frac{(L_f + L_{\hat{g}}) D^2}{k}
      +
      \frac{\sigma_* D}{\sqrt{k}}
      +
      \delta_f
      +
      \delta_{\hat{g}}
    $
           &
           &
    \cref{as:ApproximateSmoothness,as:BoundedFeasibleSet,as:ApproximatelySmoothVariance}
           &
    \cref{th:UniversalSgd-VarianceAtSolution}
    \\
    \midrule
    \multirow{2}{*}{$\UniversalFastSgd$ (\cref{alg:UniversalFastSgd})}
           &
    $\frac{L_f D^2}{k^2} + \frac{\sigma D}{\sqrt{k}} + k \delta_f$
           &
    \multirow{2}{*}{$k$}
           &
    \cref{as:ApproximateSmoothness,as:BoundedFeasibleSet,as:UniformlyBoundedVariance}
           &
    \cref{th:UniversalFastSgd}
    \\
           &
    $
      \frac{L_f D^2}{k^2}
      +
      \frac{L_{\hat{g}} D^2}{k}
      +
      \frac{\sigma_* D}{\sqrt{k}}
      +
      k \delta_f
      +
      \delta_{\hat{g}}
    $
           &
           &
    \cref{as:ApproximateSmoothness,as:BoundedFeasibleSet,as:ApproximatelySmoothVariance}
           &
    \cref{th:UniversalFastSgd-VarianceAtSolution}
    \\
    \midrule
    $\UniversalSvrg$ (\cref{alg:UniversalSvrg})
           &
    $\frac{(L_f + L_{\hat{g}}) D^2}{2^t} + \delta_f + \delta_{\hat{g}}$
           &
    $2^t + n \log t$
           &
    \cref{as:ApproximateSmoothness,as:BoundedFeasibleSet,as:ApproximatelySmoothVariance,as:ApproximatelySmoothVariance-Extra}
           &
    \cref{th:UniversalSvrg}
    \\
    \midrule
    $\UniversalFastSvrg$ (\cref{alg:UniversalFastSvrg})
           &
    $\frac{(L_f + L_{\hat{g}}) D^2}{n (t - \log \log n)} + t (\delta_f +
        \delta_{\hat{g}})$
           &
    $n t$
           &
    \cref{as:ApproximateSmoothness,as:BoundedFeasibleSet,as:ApproximatelySmoothVariance}
           &
    \cref{th:UniversalFastSvrg}
    \\
    \bottomrule
  \end{tabular}
\end{table}

  \section{Preliminaries}

\paragraph{Notation.}

We work in the space~$\RealField^d$ equipped with the standard inner product
$\InnerProduct{\cdot}{\cdot}$ and the certain Euclidean norm:
$\Norm{x} \DefinedEqual \InnerProduct{B x}{x}^{1 / 2}$,
where $B$ is a fixed positive definite matrix.
The dual norm is defined in the standard way:
$
  \DualNorm{s}
  \DefinedEqual
  \max_{\Norm{x} = 1} \InnerProduct{s}{x}
  =
  \InnerProduct{s}{B^{-1} s}.
$

For a convex function $\Map{\psi}{\RealField^d}{\RealFieldPlusInfty}$,
its (effective) domain is the following set:
$
  \EffectiveDomain \psi
  \DefinedEqual
  \SetBuilder{x \in \RealField^d}{\psi(x) < +\infty}
$.
By $\Subdifferential \psi(x)$, we denote the subdifferential of $\psi$ at a
point~$x \in \EffectiveDomain \psi$; the specific subgradients are typically
denoted by $\Gradient \psi(x)$.

A convex function~$\Map{f}{\RealField^d}{\RealField}$ is called
$(\nu, H)$-Hölder smooth for some $\nu \in \ClosedClosedInterval{0}{1}$
and $H \geq 0$ iff
$\DualNorm{\Gradient f(x) - \Gradient f(y)} \leq H \Norm{x - y}^{\nu}$
for all $x, y \in \RealField^d$ and all
$\Gradient f(x) \in \Subdifferential f(x)$,
$\Gradient f(y) \in \Subdifferential f(y)$.
Apart from the special case of $\nu = 0$, such a function~$f$ is differentiable
at every point, i.e., $\Subdifferential f(x)$ is a singleton.
A $(1, L)$-Hölder smooth function is usually called $L$-smooth.

For a convex function $\Map{\psi}{\RealField^d}{\RealFieldPlusInfty}$,
point~$x \in \RealField^d$, vector~$g \in \RealField^d$,
and coefficient~$M \geq 0$, by
$
  \ProximalMap_{\psi}(x, g, M)
  \DefinedEqual
  \argmin_{y \in \EffectiveDomain \psi} \{
    \InnerProduct{g}{y} + \psi(y) + \frac{M}{2} \Norm{y - x}^2
  \}
$,
we denote the proximal mapping.
When $M = 0$, we allow the solution to be chosen arbitrarily.

For a convex function~$\Map{f}{\RealField^d}{\RealField}$,
points $x, y \in \RealField^d$ and $\Gradient f(x) \in \Subdifferential f(x)$,
we denote the Bregman distance by
$
  \BregmanDistanceWithSubgradient{f}{\Gradient f(x)}(x, y)
  \DefinedEqual
  f(y) - f(x) - \InnerProduct{\Gradient f(x)}{y - x}
  \ (\geq 0)
$.
When the specific subgradient $\Gradient f(x)$ is clear from the context,
we use the simplified notation $\BregmanDistance{f}(x, y)$.

The positive part of~$t \in \RealField$ is denoted by
$\PositivePart{t} \DefinedEqual \max\Set{t, 0}$.
For $\tau > 0$, we also use
$\SimplifiedLog \tau \DefinedEqual \max\Set{1, \log \tau}$.

\paragraph{Problem Formulation.}

In this paper, we consider the composite optimization problem
\begin{equation}
  \label{eq:Problem}
  F^*
  \DefinedEqual
  \min_{x \in \EffectiveDomain \psi}
  \bigl[ F(x) \DefinedEqual f(x) + \psi(x) \bigr],
\end{equation}
where $\Map{f}{\RealField^d}{\RealField}$ is a convex function,
and $\Map{\psi}{\RealField^d}{\RealFieldPlusInfty}$ is a proper closed convex
function which is assumed to be sufficiently simple in the sense that the
proximal mapping $\ProximalMap_{\psi}$ can be easily computed.
We assume that this problem has a solution which we denote by~$x^*$.

To quantify the smoothness level of the objective function, we use the
following assumption:

\begin{assumption}
  \label{as:ApproximateSmoothness}
  The function~$f$ in problem~\eqref{eq:Problem} is approximately smooth:
  there exist constants $L_f, \delta_f \geq 0$ and
  $\Map{\bar{f}}{\RealField^d}{\RealField}$,
  $\Map{\bar{g}}{\RealField^d}{\RealField^d}$
  such that, for any $x, y \in \RealField^d$,
  $
    \beta_{f, \bar{f}, \bar{g}}(x, y)
    \DefinedEqual
    f(y) - \bar{f}(x) - \InnerProduct{\bar{g}(x)}{y - x}
  $
  satisfies the following inequality:
  $
    0
    \leq
    \beta_{f, \bar{f}, \bar{g}}(x, y)
    \leq
    \frac{L_f}{2} \Norm{x - y}^2 + \delta_f.
  $
\end{assumption}

\Cref{as:ApproximateSmoothness} is well-known in the literature under the
name \emph{$(\delta, L)$-oracle} and was originally introduced
in~\cite{devolder2013first}.
It covers many interesting examples.
For instance, if $f$ is $L$-smooth, then \cref{as:ApproximateSmoothness}
is satisfied with $\bar{f} = f$, $\bar{g} = \Gradient f$, $\delta_f = 0$
and $L_f = L$.
More generally, if the function $f$ is $(\nu, H_f(\nu))$-Hölder smooth, then
\cref{as:ApproximateSmoothness} is satisfied
with $\bar{f} = f$, $\bar{g} = \Gradient f$
(arbitrary selection of subgradients),
any $\delta_f > 0$ and $L_f = L(\delta_f, \nu, H_f(\nu))$, where
$
  L(\delta, \nu, H)
  \DefinedEqual
  [\frac{1 - \nu}{2 (1 + \nu) \delta}]^{\frac{1 - \nu}{1 + \nu}}
  H^{\frac{2}{1 + \nu}}
$
(see \cref{th:HolderGradientImpliesApproximateSmoothness}).
If $f$ can be uniformly approximated by an $L$-smooth function~$\phi$,
namely, $\phi(x) \leq f(x) \leq \phi(x) + \delta$, then
\cref{as:ApproximateSmoothness} is satisfied with
$\bar{f} = \phi$, $\bar{g} = \Gradient \phi$ and $\delta_f = \delta$.
If $f$ represents another auxiliary optimization problem with a strongly
concave objective, e.g., $f(x) = \max_u \Psi(x, u)$, whose
solution~$\bar{u}(x)$ can only be found with accuracy $\delta$, then $f$
satisfies \cref{as:ApproximateSmoothness} with
$\bar{f}(x) = \Psi(x, \bar{u}(x))$,
$\bar{g}(x) = \Gradient_u \Psi(x, \bar{u}(x))$ and $\delta_f = \delta$.
For more details and other interesting examples, we refer the reader
to~\cite{devolder2013first}.

In what follows, we assume that we have access to an unbiased stochastic
oracle~$\hat{g}$ for~$\bar{g}$.
Formally, this is a pair $\hat{g} = (g, \xi)$ consisting of a random
variable~$\xi$ and a mapping
$\Map{g}{\RealField^d \times \Image \xi}{\RealField^d}$
(with $\Image \xi$ being the image of~$\xi$).
When queried at a point $x$, the oracle automatically generates an independent
copy~$\xi$ of its randomness and then returns $\hat{g}_x = g(x, \xi)$
(notation: $\hat{g}_x \EqualRandom \hat{g}(x)$).
We call $g$ and~$\xi$ the function component and the random variable component
of~$\hat{g}$, respectively.
At this point, we only assume that our stochastic oracle~$\hat{g}$ is un
unbiased estimator of~$\bar{g}$, and later make various assumptions on its
variance.

Another important assumption on problem~\eqref{eq:Problem}, that we need in our
analysis, is the boundedness of the feasible set $\EffectiveDomain \psi$.

\begin{assumption}
  \label{as:BoundedFeasibleSet}
  There exists $D > 0$ such that $\Norm{x - y} \leq D$
  for any $x, y \in \EffectiveDomain \psi$.
\end{assumption}

\Cref{as:BoundedFeasibleSet} is rather standard in the literature on adaptive
methods for stochastic convex optimization
(see~\cite{levy2018online,kavis2019unixgrad,ene2021adaptive,dubois2022svrg,liu2022adaptive,rodomanov2024universal})
and can always be ensured with $D = 2 R_0$ whenever one has the knowledge of
an upper bound $R_0$ on the distance from the initial point~$x_0$ to the
solution~$x^*$.
To that end, it suffices to rewrite the problem~\eqref{eq:Problem} in
the following equivalent form:
$\min_{x \in \EffectiveDomain \psi_D} [f(x) + \psi_D(x)]$, where
$\psi_D$ is the sum of $\psi$ and the indicator function of the ball
$B_0 \DefinedEqual \SetBuilder{x \in \RealField^d}{\Norm{x - x_0} \leq R_0}$.
Note that this transformation keeps the function~$\psi_D$ reasonably simple as
its proximal mapping can be computed via that of~$\psi$ by solving a certain
one-dimensional nonlinear equation, which can be done very efficiently by
Newton's method (at no extra queries to the stochastic oracle);
in some special cases, the corresponding nonlinear equation can even be solved
analytically, e.g., when $\psi = 0$, the proximal mapping of~$\psi_D$ is
simply the projection on~$B_0$.

Throughout this paper, we refer to the constant~$D$ from
\cref{as:BoundedFeasibleSet} as the diameter of the feasible set, and
assume that its value is known to us.
This will be the only essential parameter in our methods.

  \section{Main Algorithms and Stepsize Update Rules}
\label{sec:MainAlgorithmsAndStepsizeUpdateRules}

We now present our two main algorithms for solving problem~\eqref{eq:Problem}:
$\UniversalSgd$ (\cref{alg:UniversalSgd}), and its accelerated version,
$\UniversalFastSgd$ (\cref{alg:UniversalFastSgd}).
Except the specific choice of the stepsize coefficients~$M_k$, both algorithms
are rather standard: the first one is the classical SGD method, and the second
one is the classical accelerated gradient method for stochastic
optimization~\cite{lan2012optimal}, also known as the Method of Similar
Triangles (see, e.g., Section~6.1.3 in~\cite{nesterov2018lectures}).

\begin{algorithm}[tb]
  \caption{
    $
      (\bar{x}_N, x_N, M_N)
      \EqualRandom
      \UniversalSgd_{\hat{g}, \psi}(x_0, M_0, N; D)
    $
  }
  \label{alg:UniversalSgd}
  \begin{algorithmic}[1]
    \Require
    Oracle~$\hat{g}$, comp.\ part~$\psi$,
    point~$x_0 \in \EffectiveDomain \psi$, coefficient~$M_0$,
    iteration limit~$N$, diameter~$D$.

    \State
    $g_0 \EqualRandom \hat{g}(x_0)$.
    \For{$k = 0, \ldots, N - 1$}
      \State
      $x_{k + 1} = \ProximalMap_{\psi}(x_k, g_k, M_k)$, \
      $g_{k + 1} \EqualRandom \hat{g}(x_{k + 1})$.
      \State
      $M_{k + 1} = M_+(M_k, D^2, x_k, x_{k + 1}, g_k, g_{k + 1})$
      % \begin{noindent}
      \Comment{{\scriptsize
        e.g.,
        $
          \eqcref{eq:AdaGradStepsizeUpdateRule}
          \sqrt{M_k^2 + \frac{1}{D^2} \DualNorm{g_{k + 1} - g_k}^2}
        $.
      }}
      % \end{noindent}
    \EndFor
    \State
    \Return $(\bar{x}_N, x_N, M_N)$,
    where $\bar{x}_N \DefinedEqual \frac{1}{N} \sum_{i = 1}^N x_i$.
  \end{algorithmic}
\end{algorithm}

\begin{algorithm}[tb]
  \caption{$\UniversalFastSgd_{\hat{g}, \psi}(x_0; D)$}
  \label{alg:UniversalFastSgd}
  \begin{algorithmic}[1]
    \Require
    Stochastic oracle~$\hat{g}$, composite part~$\psi$,
    point~$x_0 \in \EffectiveDomain \psi$, diameter~$D$.

    \State
    $v_0 = x_0$, $M_0 = A_0 = 0$.
    \For{$k = 0, 1, \ldots$}
      \State
      $a_{k + 1} = \frac{1}{2} (k + 1)$, $A_{k + 1} = A_k + a_{k + 1}$.
      \State
      $y_k = \frac{A_k}{A_{k + 1}} x_k + \frac{a_{k + 1}}{A_{k + 1}} v_k$, \
      $g_{y_k} \EqualRandom \hat{g}(y_k)$.
      \State
      $
        v_{k + 1}
        =
        \ProximalMap_{\psi}(v_k, g_{y_k}, \frac{M_k}{a_{k + 1}})
      $.
      \State
      $
        x_{k + 1}
        =
        \frac{A_k}{A_{k + 1}} x_k + \frac{a_{k + 1}}{A_{k + 1}} v_{k + 1}
      $, \
      $g_{x_{k + 1}} \EqualRandom \hat{g}(x_{k + 1})$.
      \State
      $
        M_{k + 1}
        \! = \!
        \frac{a_{k + 1}^2}{A_{k + 1}}
        M_+\bigl(
          \frac{A_{k + 1}}{a_{k + 1}^2} M_k,
          \frac{a_{k + 1}^2}{A_{k + 1}^2} D^2,
          y_k,
          x_{k + 1},
          g_{y_k},
          g_{x_{k + 1}}
        \bigr)
      $
      % \begin{noindent}
      \Comment{{\scriptsize
        \!\!\! e.g.,
        $
          \! \eqcref{eq:AdaGradStepsizeUpdateRule} \!
          \sqrt{
            M_k^2
            \! + \!
            \frac{a_{k + 1}^2}{D^2} \DualNorm{g_{x_{k + 1}} - g_{y_k}}^2
          }
        $.
      }}
      % \end{noindent}
    \EndFor
  \end{algorithmic}
\end{algorithm}

Both methods are expressed in terms of a certain abstract stepsize update
rule $M_+(\cdot)$ defined as follows.
Given the current stepsize coefficient $M \geq 0$, constant $\Omega > 0$
(the scaled squared diameter),
current point $x \in \EffectiveDomain \psi$ with the stochastic gradient
$\hat{g}_x \EqualRandom \hat{g}(x)$,
and the next random point~$\hat{x}_+ \in \EffectiveDomain \psi$
(possibly dependent on~$\hat{g}_x$)
with the stochastic gradient $\hat{g}_{x_+} \EqualRandom \hat{g}(\hat{x}_+)$,
the update rule computes
$\hat{M}_+ = M_+(M, \Omega, x, \hat{x}_+, \hat{g}_x, \hat{g}_{x_+})$
such that $\hat{M}_+ \geq M$ and the following inequality holds
for any $\bar{M} > c_2 L_f$:
\begin{equation}
  \label{eq:RequirementOnStepsizeUpdateRule}
  \begin{multlined}
    \Expectation[
      \hat{\Delta}(\hat{M}_+)
      +
      (\hat{M}_+ - M) \Omega
      +
      \beta_{f, \bar{f}, \bar{g}}(\hat{x}_+, x)
    ]
    \\
    \leq
    \frac{c_1}{\bar{M} - c_2 L_f}
    \Expectation[\Variance_{\hat{g}}(\hat{x}_+) + \Variance_{\hat{g}}(x)]
    +
    c_3 \delta_f
    +
    c_4 \Expectation\bigl\{
      \PositivePart{\min\Set{\hat{M}_+, \bar{M}} - M} \Omega
    \bigr\},
  \end{multlined}
\end{equation}
where
$
  \hat{\Delta}(\hat{M}_+)
  \DefinedEqual
  \beta_{f, \bar{f}, \bar{g}}(x, \hat{x}_+)
  +
  \InnerProduct{\bar{g}(x) - \hat{g}_x}{\hat{x}_+ - x}
  -
  \frac{\hat{M}_+}{2} \Norm{\hat{x}_+ - x}^2
$,
$c_1, c_2, c_3, c_4 > 0$ are some (absolute) constants,
$
  \Variance_{\hat{g}}(x)
  \DefinedEqual
  \Expectation_{\xi}[\DualNorm{g(x, \xi) - \bar{g}(x)}^2]
$
is the variance of~$\hat{g}$, and the expectation is taken w.r.t.\ the entire
randomness.

The main example is the following AdaGrad rule:
\begin{equation}
  \label{eq:AdaGradStepsizeUpdateRule}
  \boxed{
    \hat{M}_+
    =
    \sqrt{M^2 + \frac{1}{\Omega} \DualNorm{\hat{g}_{x_+} - \hat{g}_x}^2}.
  }
\end{equation}
For this rule, we have
$c_1 = \tfrac{5}{2}$, $c_2 = 4$, $c_3 = 6$, $c_4 = 2$
(see \cref{th:AdaGradRule}).
Another interesting example recently suggested in~\cite{rodomanov2024universal}
is $\hat{M}_+$ found from the equation
\begin{equation}
  \label{eq:ModifiedAdaGradRule}
  (\hat{M}_+ - M) \Omega
  =
  \PositivePart[\Big]{
    \InnerProduct{\hat{g}_{x_+} - \hat{g}_x}{\hat{x}_+ - x}
    -
    \frac{\hat{M}_+}{2} \Norm{\hat{x}_+ - x}^2
  }.
\end{equation}
For this rule, we have $c_1 = 1$, $c_2 = 2$, $c_3 = 6$, $c_4 = 2$
(see \cref{th:ModifiedAdaGradRule}).

Inequality~\eqref{eq:RequirementOnStepsizeUpdateRule} is the only property we
need from the stepsize update rule to establish all forthcoming results.
This inequality is exactly what is typically used inside the convergence proofs
for stochastic gradient methods with predefined stepsizes $M_k \equiv \bar{M}$
(in which case $M = \hat{M}_+ = \bar{M}$), where $\bar{M}$ depends on
problem-dependent constants.
The key property of AdaGrad stepsizes
(either~\eqref{eq:AdaGradStepsizeUpdateRule} or~\eqref{eq:ModifiedAdaGradRule})
is that they ensure the same inequality but now $\bar{M}$ is
the virtual stepsize existing only in the theoretical analysis.
The price for this is the extra error term
$\PositivePart{\min\Set{\hat{M}_+, \bar{M}} - M} \Omega$
appearing in the right-hand side of \cref{eq:RequirementOnStepsizeUpdateRule}.
The crucial property of this error term is that it is telescopic,
$
  \sum_{i = 0}^k \PositivePart{\min\Set{M_{i + 1}, \bar{M}} - M_i} \Omega
  =
  \PositivePart{\min\Set{M_{k + 1}, \bar{M}} - M_0} \Omega
$
(see \cref{th:TelescopingDifferencesWithMin}) and therefore its total
cumulative
impact is always bounded by the controllable constant $\bar{M} \Omega$.
Although a number of other works on theoretical analysis of AdaGrad methods
for smooth optimization use some similar ideas about the virtual stepsize
(e.g.,~\cite{levy2018online,kavis2019unixgrad,liu2022adaptive}),
this is the first time one has abstracted away all the technical details and
identified the specific inequality~\eqref{eq:RequirementOnStepsizeUpdateRule}
responsible for the universality of AdaGrad methods.

  \section{Uniformly Bounded Variance}
\label{sec:UniformlyBoundedVariance}

In this section, we assume that the variance of our stochastic oracle is
uniformly bounded.

\begin{assumption}
  \label{as:UniformlyBoundedVariance}
  For the stochastic oracle~$\hat{g}$, we have
  $
    \sigma^2
    \DefinedEqual
    \sup_{x \in \EffectiveDomain \psi} \Variance_{\hat{g}}(x)
    <
    +\infty
  $,
  where
  $
    \Variance_{\hat{g}}(x)
    \DefinedEqual
    \Expectation_{\xi}[\DualNorm{g(x, \xi) - \bar{g}(x)}^2]
  $.
\end{assumption}

Under this assumption, we can establish the following efficiency estimates for
our $\UniversalSgd$ and $\UniversalFastSgd$ methods
(the proofs are deferred to \cref{sec:UniformlyBoundedVariance-Proofs}).

\begin{restatable}{theorem}{thComplexityOfUniversalSgd}
  \label{th:UniversalSgd}
  Let \cref{alg:UniversalSgd} with $M_0 = 0$
  be applied to problem~\eqref{eq:Problem} under
  \cref{as:BoundedFeasibleSet,as:ApproximateSmoothness,as:UniformlyBoundedVariance}.
  Then, for the point $\bar{x}_N$ generated by the algorithm, we have
  \[
    \Expectation[F(\bar{x}_N)] - F^*
    \leq
    \frac{c_2 c_4 L_f D^2}{N}
    +
    2 \sigma D \sqrt{\frac{2 c_1 c_4}{N}}
    +
    c_3 \delta_f.
  \]
\end{restatable}

\begin{restatable}{theorem}{thComplexityOfUniversalFastSgd}
  \label{th:UniversalFastSgd}
  Let \cref{alg:UniversalFastSgd} be applied to
  problem~\eqref{eq:Problem} under
  \cref{as:BoundedFeasibleSet,as:ApproximateSmoothness,as:UniformlyBoundedVariance}.
  Then, for any $k \geq 1$, we have
  \[
    \Expectation[F(x_k)] - F^*
    \leq
    \frac{4 c_2 c_4 L_f D^2}{k (k + 1)}
    +
    4 \sigma D \sqrt{\frac{2 c_1 c_4}{3 k}}
    +
    \frac{c_3}{3} (k + 2) \delta_f.
  \]
\end{restatable}

We see that, in contrast to $\UniversalSgd$, the accelerated algorithm
$\UniversalFastSgd$ is not robust to the oracle's errors:
it accumulates them with time at the rate of~$\BigO(k \delta)$.
This is not surprising since the same phenomenon also occurs in the classical
accelerated gradient method, even when the oracle is deterministic and
the algorithm has the knowledge about all constants
(see~\cite{devolder2013first}).

The complexity results from \cref{th:UniversalSgd,th:UniversalFastSgd}
are similar to those from~\cite{devolder2011stochastic}.
However, it is important that our methods are adaptive and do not require
knowing the constants~$L_f$ and~$\sigma$.

In the specific case when $\delta_f = 0$, we recover the same convergence
rates as in~\cite{levy2018online,kavis2019unixgrad}, although our methods
work for the more general composite optimization problem and, in contrast
to~\cite{levy2018online}, do not require that $\Gradient f(x^*) = 0$.

  \section{Implicit Variance Reduction}
\label{sec:ImplicitVarianceReduction}

The assumption of uniformly bounded variance may not hold for some problems,
or the corresponding constant $\sigma^2$ might be quite large,
which is why there has recently been a growing interest in various alternative
variance bound assumptions~\cite{%
  moulines2011non,%
  bottou2018optimization,%
  gower2019sgd,%
  stich2019unified,%
  woodworth2021even,%
  juditsky2023sparse,%
  ilandarideva2023accelerated%
}.
One interesting option is expressing complexity bounds via the variance at the
minimizer, $\sigma_*^2 \DefinedEqual \Variance_{\hat{g}}(x^*)$,
assuming that the stochastic oracle~$\hat{g}$ satisfies some extra smoothness
conditions.
Let us show that, for our \cref{alg:UniversalSgd,alg:UniversalFastSgd},
we can also establish such bounds, moreover, this can be done
\emph{without any modifications to the algorithms}.

In this section, we study problem~\eqref{eq:Problem} under
\cref{as:BoundedFeasibleSet,as:ApproximateSmoothness} and also under the
following additional smoothness assumption on the variance:

\begin{assumption}
  \label{as:ApproximatelySmoothVariance}
  There exist $\delta_{\hat{g}}, L_{\hat{g}} \geq 0$ such that
  $
    \Variance_{\hat{g}}(x, y)
    \leq
    2 L_{\hat{g}} [\beta_{f, \bar{f}, \bar{g}}(x, y) + \delta_{\hat{g}}]
  $
  for any $x, y \in \RealField^d$, where
  $
    \Variance_{\hat{g}}(x, y)
    \DefinedEqual
    \Expectation_{\xi}[
      \DualNorm{[g(x, \xi) - g(y, \xi)] - [\bar{g}(x) - \bar{g}(y)]}^2
    ]
  $.
\end{assumption}

Note that $\Variance_{\hat{g}}(x, y)$ is the usual variance of the estimator
$g(x, \xi) - g(y, \xi)$.
Hence,
$
  \Variance_{\hat{g}}(x, y)
  \leq
  \Expectation[\DualNorm{g(x, \xi) - g(y, \xi)}^2]
$
for any $x, y$.
Furthermore, if $\hat{g}_b$ is the mini-batch version of~$\hat{g}$ of
size~$b$
(i.e., the average of $b$ i.i.d.\ samples of~$\hat{g}(x)$ at any point~$x$),
then $\Variance_{\hat{g}_b}(x, y) = \frac{1}{b} \Variance_{\hat{g}}(x, y)$
for any $x, y$.

For instance, if $f(x) = \Expectation_{\xi}[f_{\xi}(x)]$, where each
function $f_{\xi}$ is convex and $(\delta_{\xi}, L_{\xi})$-approximately
smooth with components $(\bar{f}_{\xi}, \bar{g}_{\xi})$, then,
the stochastic gradient oracle~$\hat{g}$, defined by
$g(x, \xi) \DefinedEqual \bar{g}_{\xi}(x)$ satisfies
\cref{as:ApproximatelySmoothVariance}
with $\bar{f}(x) = \Expectation_{\xi}[\bar{f}_{\xi}(x)]$,
$\bar{g}(x) = \Expectation_{\xi}[\bar{g}_{\xi}(x)]$, and
$
  \delta_{\hat{g}}
  =
  \frac{1}{L_{\max}} \Expectation_{\xi}[L_{\xi} \delta_{\xi}]
$
$(\leq \Expectation_{\xi}[\delta_{\xi}])$, $L_{\hat{g}} = L_{\max}$,
where $L_{\max} \DefinedEqual \sup_{\xi} L_{\xi}$
(see \cref{th:StandardStochasticGradientOracleIsApproximatelySmooth}).
Furthermore, if $\hat{g}_b$ is the mini-batch version of~$\hat{g}$
of size~$b$, then $\hat{g}_b$ satisfies \cref{as:ApproximatelySmoothVariance}
with the same $\delta_{\hat{g}_b} = \delta_{\hat{g}}$
but $L_{\hat{g}_b} = \frac{1}{b} L_{\hat{g}} = \frac{1}{b} L_{\max}$
which can be much smaller than $L_{\max}$ when $b$ is large enough.

Under the new assumption on the variance, $\UniversalSgd$ enjoys the following
convergence rate
(see \cref{sec:UniversalSgd-VarianceAtSolution-Proofs} for the proof).

\begin{restatable}{theorem}{thUniversalSgdVarianceAtSolution}
  \label{th:UniversalSgd-VarianceAtSolution}
  Let \cref{alg:UniversalSgd} with $M_0 = 0$
  be applied to problem~\eqref{eq:Problem} under
  \cref{as:BoundedFeasibleSet,as:ApproximateSmoothness,as:ApproximatelySmoothVariance},
  and let $\sigma_*^2 \DefinedEqual \Variance_{\hat{g}}(x^*)$.
  Then, for the point $\bar{x}_N$ produced by the method, we have
  \[
    \Expectation[F(\bar{x}_N)] - F^*
    \leq
    \frac{c_4 (c_2 L_f + 12 c_1 L_{\hat{g}}) D^2}{N}
    +
    2 \sigma_* D \sqrt{\frac{6 c_1 c_4}{N}}
    +
    c_3 \delta_f + \frac{4}{3} \delta_{\hat{g}}.
  \]
\end{restatable}

Comparing the above result with \cref{th:UniversalSgd}, we see that
we have essentially replaced the uniform bound~$\sigma$ with the more refined
one~$\sigma_*$ at the cost of replacing $L_f$ with $L_f + L_{\hat{g}}$ and
$\delta_f$ with $\delta_f + \delta_{\hat{g}}$.
This corresponds to classical results on the usual SGD for which we know
all problem dependent-constants.
However, our method is universal and works automatically under both
assumptions from the previous section and the current one, and therefore
enjoys the best among the rates given by
\cref{th:UniversalSgd,th:UniversalSgd-VarianceAtSolution}.

For the accelerated algorithm $\UniversalFastSgd$, we have the following result
(whose proof is located in
  \cref{sec:UniversalFastSgd-VarianceAtSolution-Proofs}).

\begin{restatable}{theorem}{thUniversalFastSgdVarianceAtSolution}
  \label{th:UniversalFastSgd-VarianceAtSolution}
  Let \cref{alg:UniversalFastSgd} be applied to
  problem~\eqref{eq:Problem} under
  \cref{as:BoundedFeasibleSet,as:ApproximateSmoothness,as:ApproximatelySmoothVariance},
  and let $\sigma_*^2 \DefinedEqual \Variance_{\hat{g}}(x^*)$.
  Then, for any $k \geq 1$, we have
  \[
    \Expectation[F(x_k)] - F^*
    \leq
    \frac{4 c_2 c_4 L_f D^2}{k (k + 1)}
    +
    \frac{24 c_1 c_4 L_{\hat{g}} D^2}{k + 1}
    +
    4 \sigma_* D \sqrt{\frac{2 c_1 c_4}{k}}
    +
    \frac{c_3}{3} (k + 2) \delta_f + \frac{4}{3} \delta_{\hat{g}}.
  \]
\end{restatable}

Comparing our previous complexity bound for $\UniversalFastSgd$ under the
assumption on uniformly bounded variance (\cref{th:UniversalFastSgd})
with the bound from \cref{th:UniversalFastSgd-VarianceAtSolution},
we see that, instead of simply replacing $\sigma$ with $\sigma_*$,
$L_f$ with $L_f + L_{\hat{g}}$ and $\delta_f$ with
$\delta_f + \delta_{\hat{g}}$, which was the case for the basic method,
the situation is now not that simple.
Specifically, the $L_f$ and $L_{\hat{g}}$ terms now converge at different
rates: $\BigO(\frac{1}{k^2})$ and $\BigO(\frac{1}{k})$, respectively.
While this may seem strange at first, this behavior is actually unavoidable,
at least in the case when $\delta_f = \delta_{\hat{g}} = 0$
(see, e.g., Section~E in~\cite{woodworth2021even}).
For the case when $\delta_f = \delta_{\hat{g}} = 0$, the complexity result
from \cref{th:UniversalFastSgd-VarianceAtSolution} is similar to the results
for the Accelerated SGD algorithm from~\cite{woodworth2021even}.
However, the latter paper studies a specific setting where
$f(x) = \Expectation[f_{\xi}(x)]$, where each component $f_{\xi}$
is $L_{\max}$-smooth and then assumes that $f$ is also $L_{\max}$-smooth,
instead of working with the constant $L_f$ which can be much smaller than
$L_{\max}$.
A similar separation of the constants $L_f$ and $L_{\hat{g}}$, which we do,
was recently considered in~\cite{ilandarideva2023accelerated},
where the authors obtained some similar rates to our
\cref{th:UniversalFastSgd-VarianceAtSolution}.
However, it is important that, unlike the algorithms considered
in~\cite{woodworth2021even,ilandarideva2023accelerated},
our $\UniversalFastSgd$ is universal and does not require knowing any
problem-dependent constants except~$D$.
Furthermore, our results are more general because we allow the oracle to be
inexact\footnote{%
  One should mention that \cite{ilandarideva2023accelerated} studies a slightly
  more general assumption than
  $\Variance_{\hat{g}}(x, y) \leq 2 L_{\hat{g}} \BregmanDistance{f}(x, y)$
  (which is our \cref{as:ApproximatelySmoothVariance} for the exact oracle
    with $\delta = 0$).
  Specifically, their assumption is
  $
    \Variance_{\hat{g}}(x)
    \leq
    2 L_{\hat{g}} \BregmanDistance{f}(x^*, x)
    +
    \sigma_*^2
  $.
  While the former assumption implies the latter, they do not seem to be
  completely equivalent.
  Under the latter assumption, the authors
  in~\cite{ilandarideva2023accelerated}
  proved the same convergence rate (up to absolute constants) as in
  our \cref{th:UniversalFastSgd-VarianceAtSolution}
  with $\delta_f = \delta_{\hat{g}} = 0$.
  However, instead of using the Similar Triangles version of the Fast Gradient
  Method, they had to consider an alternative algorithm.
  We expect that this alternative algorithm can also be made universal
  using the ideas from our paper, and leave the details to a future work.
}.

  \section{Explicit Variance Reduction with SVRG}
\label{sec:ExplicitVarianceReductionWithSvrg}

Let us now show that we can also incorporate explicit SVRG-type variance
reduction into our methods.
In this section, we consider problem~\eqref{eq:Problem} under
\cref{%
  as:BoundedFeasibleSet,%
  as:ApproximateSmoothness,%
  as:ApproximatelySmoothVariance%
}.
All the proofs are deferred to
\cref{sec:ExplicitVarianceReductionWithSvrg-Proofs}.

In addition to the stochastic oracle~$\hat{g}$, we now assume that we can
also compute the (approximate) full-gradient oracle~$\bar{g}$.
This allows us to define the following auxiliary \emph{SVRG oracle} induced by
$\hat{g}$ with center $\tilde{x} \in \RealField^d$
(notation $\hat{G} = \SvrgOracle_{\hat{g}, \bar{g}}(\tilde{x})$)
as the oracle with the same random variable component~$\xi$ as~$\hat{g}$
and the function component given by
$G(x, \xi) = g(x, \xi) - g(\tilde{x}, \xi) + \bar{g}(\tilde{x})$.

Our $\UniversalSvrg$ method is presented in \cref{alg:UniversalSvrg}.
This is the classical epoch-based SVRG algorithm which can be seen as
the adaptive version of the SVRG++ method from~\cite{allen2016improved}.
A similar scheme was suggested in~\cite{dubois2022svrg}, however, instead
of accumulating gradient differences as
in~\eqref{eq:AdaGradStepsizeUpdateRule},
their method accumulates gradients and therefore does not work without the
additional assumption that $\Gradient f(x^*) = 0$
(which may not hold for constrained optimization).

\begin{algorithm}[tb]
  \caption{$\UniversalSvrg_{\hat{g}, \bar{g}, \psi}(x_0; D)$}
  \label{alg:UniversalSvrg}
  \begin{algorithmic}[1]
    \Require
    Oracles~$\hat{g}$, $\bar{g}$, composite part~$\psi$,
    point $x_0 \in \EffectiveDomain \psi$, diameter~$D$.

    \State
    $\tilde{x}_0 = x_0$, $M_0 = 0$.
    \For{$t = 0, 1, \ldots$}
      \State
      $
        (\tilde{x}_{t + 1}, x_{t + 1}, M_{t + 1})
        \EqualRandom
        \UniversalSgd_{\hat{G}_t, \psi}(x_t, M_t, 2^{t + 1}; D)
      $
      with $\hat{G}_t = \SvrgOracle_{\hat{g}, \bar{g}}(\tilde{x}_t)$.
    \EndFor
  \end{algorithmic}
\end{algorithm}

Let us now present the complexity guarantees.
To do so, we first need to introduce, one more assumption we need in our
analysis.

\begin{assumption}
  \label{as:ApproximatelySmoothVariance-Extra}
  The variance of~$\hat{g}$ satisfies
  $
    \Variance_{\hat{g}}(x, y)
    \leq
    4 L_{\hat{g}} [
      \BregmanDistanceWithSubgradient{f}{\Gradient f(x)}(x, y)
      +
      2 \delta_{\hat{g}}
    ]
  $
  for any $x, y \in \RealField^d$ and
  any $\Gradient f(x) \in \Subdifferential f(x)$.
\end{assumption}

\Cref{as:ApproximatelySmoothVariance-Extra} is very similar to
\cref{as:ApproximatelySmoothVariance}.
The only difference between them is that the former contains the standard
Bregman distance in the right-hand side, while the latter contains its
approximation $\beta_{f, \bar{f}, \bar{g}}(x, y)$ involving the approximate
function value~$\bar{f}(x)$ and the approximate gradient~$\bar{g}(x)$.
Nevertheless, both assumptions are actually satisfied for the main examples
we discussed after introducing \cref{as:ApproximatelySmoothVariance}
(see \cref{th:StandardStochasticGradientOracleIsApproximatelySmooth}).

\begin{restatable}{theorem}{thUniversalSvrg}
  \label{th:UniversalSvrg}
  Let $\UniversalSvrg$ (as defined by \cref{alg:UniversalSvrg})
  be applied to problem~\eqref{eq:Problem} under
  \cref{%
    as:BoundedFeasibleSet,%
    as:ApproximateSmoothness,%
    as:ApproximatelySmoothVariance,%
    as:ApproximatelySmoothVariance-Extra%
  }.
  Then, for any $t \geq 1$ and $\bar{c}_3 \DefinedEqual \max\Set{c_3, 1}$,
  we have
  \[
    \Expectation[F(\tilde{x}_t)] - F^*
    \leq
    \frac{[(c_2 c_4 + 1) L_f + 48 c_1 c_4 L_{\hat{g}}] D^2}{2^t}
    +
    2 \bar{c}_3 \delta_f
    +
    \frac{8}{3} \delta_{\hat{g}}.
  \]
  To construct~$\tilde{x}_t$, the algorithm needs to make $\BigO(2^t)$ queries
  to~$\hat{g}$ and $\BigO(t)$ queries to~$\bar{g}$.
\end{restatable}

We now present an accelerated version of $\UniversalSvrg$,
see \cref{alg:UniversalFastSvrg}.
As $\UniversalSvrg$, this method is also epoch-based, and its epoch is very
similar to $\UniversalFastSgd$ (\cref{alg:UniversalFastSvrg}) in the sense that
it also iterates similar-triangle steps.
However, the triangles in $\UniversalTriangleSvrgEpoch$ are of the form
$(\tilde{x}, v_k, v_{k + 1})$, i.e., they always share the common
vertex~$\tilde{x}$, in contrast to the triangles $(x_k, v_k, v_{k + 1})$
in $\UniversalFastSgd$ (in $\UniversalTriangleSvrgEpoch$, the role of the
  average points~$y_k$ is played by~$x_k$).
We note that our $\UniversalFastSvrg$ is essentially the primal version of the
VRADA method from~\cite{song2020variance}, but equipped with AdaGrad stepsizes.
Alternative accelerated SVRG schemes with AdaGrad
stepsizes~\eqref{eq:AdaGradStepsizeUpdateRule} were recently proposed
in~\cite{liu2022adaptive}; however, they seem to be much more complicated.

The special choice of the initial reference point~$\tilde{x}_0$
at \cref{step:FastSvrg-ChooseInitialPoints} is rather standard
and motivated by the desire to keep the initial function residual appropriately
bounded: $F(\tilde{x}_0) - F^* \leq \frac{1}{2} L_f D^2 + \delta_f$;
the simplest way to achieve this is to make the full gradient step from any
feasible point (see \cref{th:FullGradientStepGivesGoodInitialPoint}).

\begin{figure}[tb]
  \captionsetup{hypcap=false}
  \begin{algorithm}[H]
    \caption{$\UniversalFastSvrg_{\hat{g}, \bar{g}, \psi}(x_0, N; D)$}
    \label{alg:UniversalFastSvrg}
    \begin{algorithmic}[1]
      \Require
      Oracles~$\hat{g}$, $\bar{g}$, composite part~$\psi$,
      point~$x_0 \in \EffectiveDomain \psi$, epoch length~$N$, diameter~$D$.

      \State
      \label{step:FastSvrg-ChooseInitialPoints}%
      $\tilde{x}_0 = \ProximalMap_{\psi}(x_0, \bar{g}(x_0), 0)$,
      $v_0 = x_0$, $M_0 = 0$, $A_0 = \frac{1}{N}$.
      \For{$t = 0, 1, \ldots$}
        \State
        $a_{t + 1} = \sqrt{A_t}$, \ $A_{t + 1} = A_t + a_{t + 1}$.
        \State
        $
          (\tilde{x}_{t + 1}, v_{t + 1}, M_{t + 1})
          \EqualRandom
          \UniversalTriangleSvrgEpoch_{\hat{g}, \bar{g}, \psi}(
            \tilde{x}_t, v_t, M_t, A_t, a_{t + 1}, N; D
          )
        $.
      \EndFor
    \end{algorithmic}
  \end{algorithm}
  \vspace{-2em}
  \begin{algorithm}[H]
    \caption{
      $
        (\tilde{x}_+, v_+, M_+)
        \EqualRandom
        \UniversalTriangleSvrgEpoch_{\hat{g}, \bar{g}, \psi}(
          \tilde{x}, v_0, M_0, A, a, N; D
        )
      $
    }
    \label{alg:UniversalTriangleSvrgEpoch}
    \begin{algorithmic}[1]
      \Require
      Oracles~$\hat{g}$, $\bar{g}$, comp.\ part~$\psi$,
      points~$\tilde{x}, v_0$, coefficients~$M_0, A, a$,
      epoch length~$N$, diameter~$D$.

      \State
      $A_+ = A + a$, $x_0 = \frac{A}{A_+} \tilde{x} + \frac{a}{A_+} v_0$,
      $\hat{G} = \SvrgOracle_{\hat{g}, \bar{g}}(\tilde{x})$,
      $G_{x_0} \EqualRandom \hat{G}(x_0)$.
      \For{$k = 0, \ldots, N - 1$}
        \State
        $v_{k + 1} = \ProximalMap_{\psi}(v_k, G_{x_k}, \frac{M_k}{a})$.
        \State
        $x_{k + 1} = \frac{A}{A_+} \tilde{x} + \frac{a}{A_+} v_{k + 1}$, \
        $G_{x_{k + 1}} \EqualRandom \hat{G}(x_{k + 1})$.
        \State
        $
          M_{k + 1}
          =
          \frac{a^2}{A_+}
          M_+\bigl(
            \frac{A_+}{a^2} M_k,
            \frac{a^2}{A_+^2} D^2,
            x_k,
            x_{k + 1},
            G_{x_k},
            G_{x_{k + 1}}
          \bigr)
        $
        % \begin{noindent}
        \Comment{{\scriptsize
          e.g.,
          $
            \eqcref{eq:AdaGradStepsizeUpdateRule}
            \sqrt{
              M_k^2
              +
              \frac{a^2}{D^2} \DualNorm{G_{x_{k + 1}} - G_{x_k}}^2
            }
          $.
        }}
        % \end{noindent}
      \EndFor
      \State
      \Return $(\bar{x}_N, v_N, M_N)$,
      where $\bar{x}_N \DefinedEqual \frac{1}{N} \sum_{k = 1}^N x_k$.
    \end{algorithmic}
  \end{algorithm}
\end{figure}

\begin{restatable}{theorem}{ThUniversalFastSvrg}
  \label{th:UniversalFastSvrg}
  Let $\UniversalFastSvrg$ (\cref{alg:UniversalFastSvrg}) be applied to
  problem~\eqref{eq:Problem} under \cref{%
    as:BoundedFeasibleSet,%
    as:ApproximateSmoothness,%
    as:ApproximatelySmoothVariance%
  }, and let $N \geq 9$.
  Then, for any
  $t \geq t_0 \DefinedEqual \Ceil{\log_2 \log_3 N} - 1 \ (\geq 0)$,
  it holds that
  \[
    \Expectation[F(\tilde{x}_t)] - F^*
    \leq
    \frac{
      9 [(c_2 c_4 + \frac{1}{2}) L_f + 6 c_1 c_4 L_{\hat{g}}] D^2
    }{
      N (t - t_0 + 1)^2
    }
    +
    (c_3 t + 1) \delta_f + \frac{5}{3} t \delta_{\hat{g}}.
  \]
  To construct~$\tilde{x}_t$, the algorithm needs to make $\BigO(N t)$ queries
  to~$\hat{g}$ and $\BigO(t)$ queries to~$\bar{g}$.
  Assuming that the complexity of querying~$\bar{g}$ is $n$ times bigger
  than that of querying~$\hat{g}$ and choosing $N = \BigTheta(n)$,
  we get the total stochastic-oracle complexity of $\BigO(n t)$.
\end{restatable}

Note that \cref{th:UniversalFastSvrg}, unlike \cref{th:UniversalSvrg},
does not require the extra \cref{as:ApproximatelySmoothVariance-Extra}.
This suggests that \cref{as:ApproximatelySmoothVariance-Extra} might be
somewhat artificial and could potentially be removed from
\cref{th:UniversalSvrg} as well.
However, we do not know how to do it, even in the simplest case when
$\delta_f = \delta_{\hat{g}} = 0$ and the algorithm has the knowledge of the
constants $L_f$ and $L_{\hat{g}}$ from
\cref{as:ApproximateSmoothness,as:ApproximatelySmoothVariance}.

  \section{Application to Hölder Smooth Problems}
\label{sec:ApplicationToHolderSmoothProblems}

To illustrate how powerful our results are, let us quickly consider the
specific example of solving the stochastic optimization problem with Hölder
smooth components.

\begin{example}
  \label{ex:ProblemWithHolderSmoothComponents}
  Suppose that the function~$f$ in problem~\eqref{eq:Problem} is the expectation
  of other functions, $f(x) = \Expectation_{\xi}[f_{\xi}(x)]$,
  where each function $f_{\xi}$ is convex and $(\nu, H_{\xi}(\nu))$-Hölder
  smooth.
  Consider the standard mini-batch stochastic gradient oracle $\hat{g}_b$ of
  size~$b$, defined by
  $g_b(x, \xi_{[b]}) = \frac{1}{b} \sum_{j = 1}^b \Gradient f_{\xi_j}(x)$,
  where $\xi_{[b]} \DefinedEqual (\xi_1, \ldots, \xi_b)$
  with $b$ i.i.d.\ copies of~$\xi$,
  and $\Gradient f_{\xi}(x) \in \Subdifferential f_{\xi}(x)$ is an arbitrary
  selection of subgradients for each~$\xi$.
  We define $H_f(\nu)$ as the Hölder constant for the function~$f$ and
  $H_{\max}(\nu) \DefinedEqual \sup_{\xi} H_{\xi}(\nu)$
  as the worst among Hölder constants for each $f_{\xi}$.
  Note that we always have $H_f(\nu) \leq \Expectation_{\xi}[H_{\xi}(\nu)]$
  but $H_f(\nu)$ can, in principle, be much smaller than the right-hand side.
  Also, define
  $
    \sigma^2
    \DefinedEqual
    \sup_{x \in \EffectiveDomain \psi} \Variance_{\hat{g}_1}(x)
    \equiv
    \sup_{x \in \EffectiveDomain \psi}
    \Expectation_{\xi}[\DualNorm{\Gradient f_{\xi}(x) - \Gradient f(x)}^2]
  $
  and
  $
    \sigma_*^2
    \DefinedEqual
    \Variance_{\hat{g}_1}(x^*)
    \equiv
    \Expectation_{\xi}[\DualNorm{\Gradient f_{\xi}(x^*) - \Gradient f(x^*)}^2]
  $.
  We assume that the computation of $\hat{g}_b$ can be parallelized
  and the computation of $\Gradient f$ is $n_b$ times more expensive than that
  of $\hat{g}_b$.
\end{example}

To solve the above problem, we can apply any of the methods we presented
before.
The resulting oracle complexities (in terms of the BigO-notation) are
summarized in \cref{tab:CorollariesForHolderSmoothProblem};
the precise statements the corresponding results and their proofs
are deferred to \cref{sec:ApplicationToHolderSmoothProblems-Proofs}.

\begin{table}
  \crefname{algorithm}{Alg.}{Algs.}
  \crefname{corollary}{Cor.}{Cors.}
  \renewcommand{\crefpairconjunction}{, }
  \centering
  \captionsetup{font=small}
  \caption{%
    Corollaries of our results for the case when problem~\eqref{eq:Problem} has
    Hölder smooth components, as defined in
    \cref{ex:ProblemWithHolderSmoothComponents}.
    ``SO complexity'' is the stochastic-oracle complexity for
    reaching accuracy~$\epsilon$ in terms of the expected function residual,
    defined as in \cref{tab:SummaryOfMainResults} but with
    $\hat{g} = \hat{g}_b$, $\bar{g} = \Gradient f$, $n = n_b$.
  }
  \label{tab:CorollariesForHolderSmoothProblem}
  \scriptsize
  \begin{tabular}{lccc}
    \toprule
    Method & SO complexity & Reference \\
    \midrule
    $\UniversalSgd$ (\cref{alg:UniversalSgd})
           &
    $
      \bigl( \frac{H_f(\nu)}{\epsilon} \bigr)^{\frac{2}{1 + \nu}} D^2
      +
      \frac{1}{b}
      \min\Set[\big]{
        % \begin{noindent}
        \frac{\sigma^2 D^2}{\epsilon^2},
        \bigl( \frac{H_{\max}(\nu)}{\epsilon} \bigr)^{\frac{2}{1 + \nu}} D^2
        +
        \frac{\sigma_*^2 D^2}{\epsilon^2}
        % \end{noindent}
      }
    $
           &
    \cref{th:UniversalSgd-Holder,th:UniversalSgd-VarianceAtSolution-Holder}
    \\
    \midrule
    $\UniversalFastSgd$ (\cref{alg:UniversalFastSgd})
           &
    $
      \bigl( \frac{H_f(\nu) D^{1 + \nu}}{\epsilon} \bigr)^{\frac{2}{1 + 3 \nu}}
      +
      \frac{1}{b}
      \min\Set[\big]{
        % \begin{noindent}
        \frac{\sigma^2 D^2}{\epsilon^2},
        \bigl( \frac{H_{\max}(\nu)}{\epsilon} \bigr)^{\frac{2}{1 + \nu}} D^2
        +
        \frac{\sigma_*^2 D^2}{\epsilon^2}
        % \end{noindent}
      }
    $
           &
    \cref{th:UniversalFastSgd-Holder,th:UniversalFastSgd-VarianceAtSolution-Holder}
    \\
    \midrule
    $\UniversalSvrg$ (\cref{alg:UniversalSvrg})
           &
    $
      \bigl[
        N_{\nu}(\epsilon)
        \DefinedEqual
        \bigl( \frac{H_f(\nu)}{\epsilon} \bigr)^{\frac{2}{1 + \nu}} D^2
        +
        \frac{1}{b}
        \bigl( \frac{H_{\max}(\nu)}{\epsilon} \bigr)^{\frac{2}{1 + \nu}}
        D^2
      \bigr]
      +
      n_b \SimplifiedLog N_{\nu}(\epsilon)
    $
           &
    \cref{th:UniversalSvrg-Holder}
    \\
    \midrule
    $\UniversalFastSvrg$ (\cref{alg:UniversalFastSvrg})
           &
    $
      [
        \frac{n_b^{\nu} H_f(\nu) D^{1 + \nu}}{\epsilon}
      ]^{\frac{2}{1 + 3 \nu}}
      +
      [
        \frac{n_b^{\nu} H_{\max}(\nu) D^{1 + \nu}}{b^{(1 + \nu) / 2} \epsilon}
      ]^{\frac{2}{1 + 3 \nu}}
      +
      n_b \log \log n_b
    $
           &
    \cref{th:UniversalFastSvrg-Holder}
    \\
    \bottomrule
  \end{tabular}
\end{table}

Note that our problem is characterized by a large number of parameters,
$\nu$, $H_f(\nu)$, $H_{\max}(\nu)$, $\sigma$, $\sigma_*$.
For each combination of these parameters, we get a certain complexity guarantee
for each of our methods, and it is impossible to say in advance which
combination results in the smaller complexity bound.
However, it is not important for our methods since none of them needs to know
any of these constants to ensure the corresponding bound.
This means that our algorithms are \emph{universal}: they automatically
figure out the best problem class for a specific problem given to them.

  \section{Experiments}
\label{sec:Experiments}

Let us illustrate the performance of our methods in preliminary numerical
experiments on solving the following test problem:
\begin{equation}
  \label{eq:PolyhedronFeasibilityProblem}
  f^*
  \DefinedEqual
  \min_{\Norm{x} \leq R} \Bigl\{
    f(x)
    \DefinedEqual
    \frac{1}{n} \sum_{i = 1}^n \PositivePart{\InnerProduct{a_i}{x} - b_i}^q
  \Bigr\},
\end{equation}
where $a_i, b_i \in \RealField^d$, $q \in \ClosedClosedInterval{1}{2}$
and $R > 0$.

This problem covers several interesting applications.
Indeed, if $q = 2$, we get the classical Least squares problem.
If $q = 1$, this is the well-known Support-Vector Machines (SVM) problem.
In both cases, the ball-constraint $\Norm{x} \leq R$ acts as a regularizer,
and problem~\eqref{eq:PolyhedronFeasibilityProblem} is, in fact, equivalent
to $\min_{x \in \RealField^d} [f(x) + \frac{\mu}{2} \Norm{x}^2]$
for a certain $\mu \geq 0$
(this follows, e.g., from the KKT optimality conditions)
such that $\mu$ decreases when $R$ increases.

Another interesting application of~\eqref{eq:PolyhedronFeasibilityProblem},
which we consider in this section, is the
\emph{polyhedron feasibility problem}:
find a point $x^* \in \RealField^d$, $\Norm{x^*} \leq R$, inside the polyhedron
$P = \SetBuilder{x}{\InnerProduct{a_i}{x} \leq b_i, \ i = 1, \ldots, n}$.
Such a point exists iff $f^* = 0$.
Note that \eqref{eq:PolyhedronFeasibilityProblem} is a problem with
Hölder smooth components of degree $\nu = q - 1$.
By varying $q$ in~\eqref{eq:PolyhedronFeasibilityProblem},
we can therefore check the adaptivity of different methods to the unknown to
them Hölder characteristics of the objective function.

The data for our problem is generated randomly.
First, we generate $x^*$ uniformly from the sphere of radius~$0.95 R$ centered
at the origin.
Then, we generate i.i.d.\ vectors $a_i$ with components uniformly distributed
on $\ClosedClosedInterval{-1}{1}$.
We then make sure that $\InnerProduct{a_n}{x^*} < 0$ by inverting the sign
of~$a_n$ if necessary.
Next, we generate positive reals $s_i$ uniformly in
$\ClosedClosedInterval{0}{-0.1 c_{\min}}$,
where $c_{\min} \DefinedEqual \min_i \InnerProduct{a_i}{x^*} < 0$,
and set $b_i = \InnerProduct{a_i}{x^*} + s_i$.
By construction, $x^*$ is a solution of our problem with $f^* = 0$,
and the origin $x_0 = 0$ lies outside the polyhedron since there exists~$j$
(corresponding to $c_{\min}$) such that
$b_j = c_{\min} + s_j \leq 0.9 c_{\min} < 0$.

\begin{figure}[tb]
  \centering
  \includegraphics[width=\textwidth]{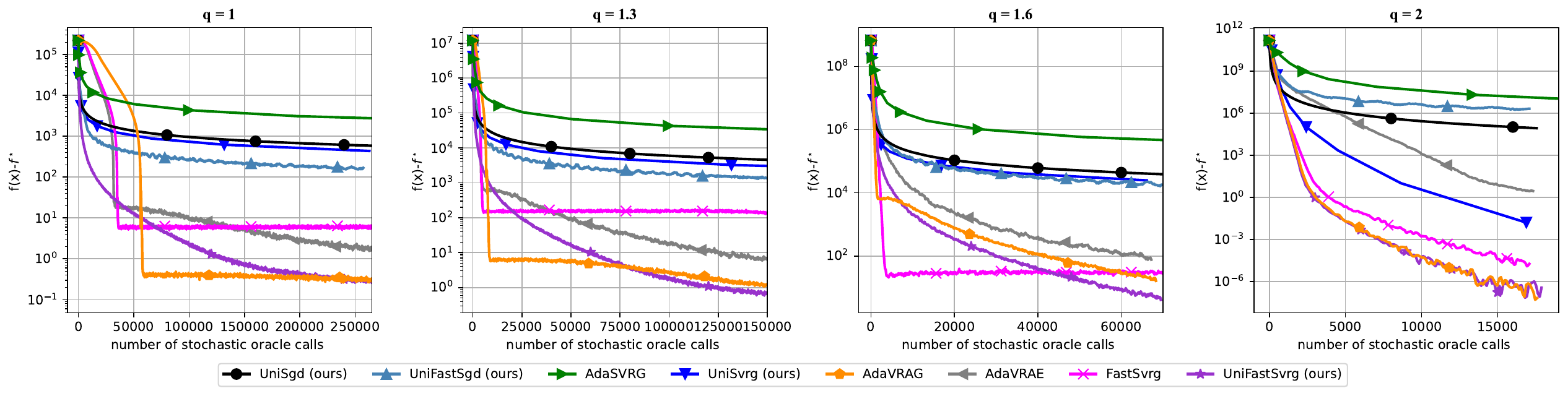}
  \caption{
    Comparison of different methods on the polyhedron feasibility problem.
  }
  \label{fig:Polyhedron}
\end{figure}

\begin{figure}[tb]
  \centering
  \includegraphics[width=\textwidth]{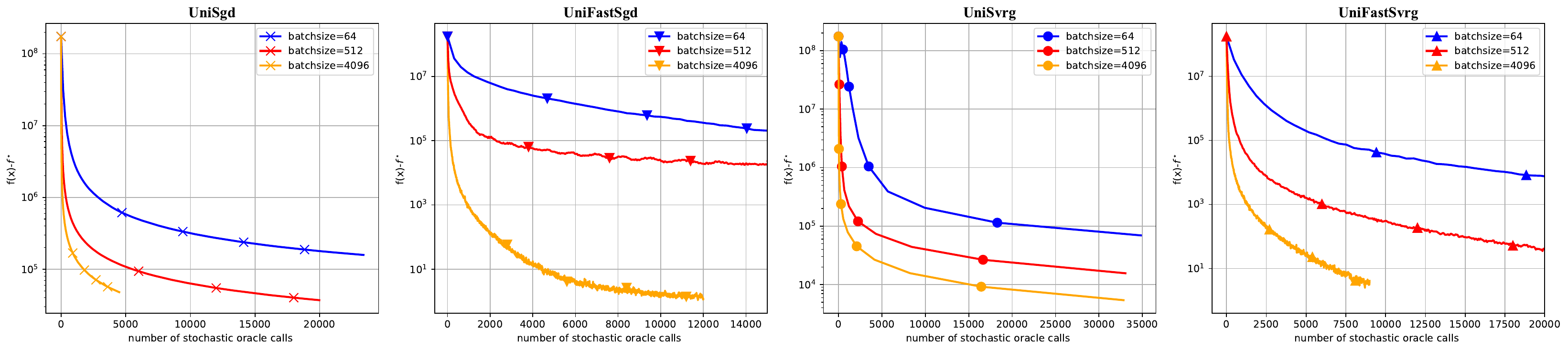}
  \caption{
    Impact of mini-batch size on the performance of our methods on the
    polyhedron feasibility problem.
  }
  \label{fig:mini-batchPolyhedron}
\end{figure}

We compare $\UniversalSvrg$ (\cref{alg:UniversalSvrg}) against
AdaSVRG~\cite{dubois2022svrg} (with parameters $K = 3$ and $\eta = D = 2 R$).
We next compare $\UniversalFastSvrg$ (\cref{alg:UniversalFastSvrg})
against AdaVRAE and AdaVRAG~\cite{liu2022adaptive}.
We also compare it with the FastSvrg method with constant
stepsize, which is the primal version of the VRADA method
from~\cite{song2020variance};
the stepsize is selected by doing a grid search over
$\SetBuilder{10^j}{j = -3, \ldots, 4}$ and choosing the best value in the
sense that the algorithm is neither too slow nor has a large error.
We report $\UniversalSgd$ (\cref{alg:UniversalSgd})
and $\UniversalFastSgd$ (\cref{alg:UniversalFastSgd}) together with these
methods.
For $\UniversalFastSvrg$, contrary to the theoretical recommendation of
choosing~$\tilde{x}_0$ as the result of the full gradient step, we found it
slightly more useful to simply set $\tilde{x}_0 = x_0$.
For all our methods, we use the AdaGrad
stepsize~\eqref{eq:AdaGradStepsizeUpdateRule};
the other stepsize~\eqref{eq:ModifiedAdaGradRule} works very similarly
(see \cref{sec:AdditionalExperiments:ComparisonBetweenStepsizeUpdateRules}
  for a detailed comparison).
For all methods, we use the standard mini-batch stochastic gradient oracle
of size $b = 256$.

The results are shown in \cref{fig:Polyhedron}, where we fix
$n = 10^4$, $d = 10^3$, $R = 10^6$ and consider different values of
$q \in \Set{1, 1.3, 1.6, 2}$.
We plot the total number of stochastic oracle calls against the function
residual.
We treat one mini-batch oracle computation as one stochastic oracle call.
If we compute the full gradient, we count this as $n / b$ stochastic oracle
calls where $n$ is the total number of samples and $b$ denotes the mini-batch
size.

We see that, except the AdaSVRG method, all SVRG algorithms typically
converge much faster than the usual SGD methods without explicit variance
reduction, at least after a few computations of the full gradient.
Among the non-accelerated SVRG methods, $\UniversalSvrg$ converges consistently
faster than AdaSVRG, while $\UniversalFastSvrg$ performs the best across the
accelerated ones.
Note that FastSvrg with constant stepsize is not converging when the problem
is not Lipschitz smooth ($q < 2$), in contrast to our universal methods.

In \cref{fig:mini-batchPolyhedron}, we also illustrate the impact of the
mini-batch size~$b$ on the convergence of our methods.
We consider the same values of $n$, $d$, $R$ as before and fix $q = 1.5$.
As we can see, in the idealized situation, when one can implement the
mini-batch oracle computations by a perfect parallelism, there is a significant
speedup in convergence when increasing the mini-batch size, exactly as
predicted by our theory.

For additional experiments, including the discussion of the implicit variance
reduction properties, we refer the reader to
\cref{sec:AdditionalExperiments}.

  \section{Conclusions}

In this paper, we showed that AdaGrad stepsizes can be applied, in a
unified manner, in a large variety of situations, leading to universal
methods suitable for multiple problem classes at the same time.
Note that this does not come for free.
We still need to know one parameter, the diameter~$D$ of the feasible set.
While it is not necessary to know this parameter precisely, the cost of
underestimating or overestimating it, can be high
(all complexity bounds would be multiplied by the ratio between our guess and
  the true~$D$).
At the same time, there already exist some parameter-free methods
which are based on AdaGrad and aim to solve precisely this problem~\cite{%
  carmon2022making,%
  defazio2023learning,%
  ivgi2023dog,%
  khaled2023dowg,%
  mishchenko2023prodigy%
}.
It is therefore interesting to consider extensions of our results to these more
advanced algorithms.
Another interesting direction is, of course, nonconvex problems.

  \printbibliography

  \newpage
  \appendix
  \section{General Auxiliary Results}

\subsection{Approximately Smooth Functions}

\begin{theorem}[Lemma~2 in~\cite{nesterov2015universal}]
  \label{th:HolderGradientImpliesApproximateSmoothness}
  Let $\Map{f}{\RealField^d}{\RealField}$ be a convex $(\nu, H)$-Hölder smooth
  function for some $\nu \in \ClosedClosedInterval{0}{1}$ and $H \geq 0$.
  Then, for any $\delta > 0$, any $x, y \in \RealField^d$ and
  any $\Gradient f(x) \in \Subdifferential f(x)$, it holds that
  $
    \BregmanDistanceWithSubgradient{f}{\Gradient f(x)}(x, y)
    \leq
    \frac{L}{2} \Norm{x - y}^2 + \delta
  $
  with
  $
    L
    =
    [\frac{1 - \nu}{2 (1 + \nu) \delta}]^{\frac{1 - \nu}{1 + \nu}}
    H^{\frac{2}{1 + \nu}}
  $
  (with the convention that $0^0 = 1$).
\end{theorem}

\begin{theorem}
  \label{th:InequalitiesForApproximatelySmoothFunctions}
  Let $\Map{f}{\RealField^d}{\RealField}$ be a $(\delta, L)$-approximately
  smooth convex function with components $(\bar{f}, \bar{g})$, i.e.,
  for any $x, y \in \RealField^d$ and
  $
    \beta_{f, \bar{f}, \bar{g}}(x, y)
    \DefinedEqual
    f(y) - \bar{f}(x) - \InnerProduct{\bar{g}(x)}{y - x}
  $,
  we have
  \begin{equation}
    \label{eq:ApproximateSmoothness-UpperBoundOnBregmanDistance}
    0
    \leq
    \beta_{f, \bar{f}, \bar{g}}(x, y)
    \leq
    \frac{L}{2} \Norm{x - y}^2 + \delta.
  \end{equation}
  Then, for any $x, y \in \RealField^d$ and any
  $\Gradient f(x) \in \Subdifferential f(x)$, the following inequalities hold:
  \begin{gather}
    \label{eq:ApproximateSmoothness-RelationBetweenFunctionValues}
    \bar{f}(x) \leq f(x) \leq \bar{f}(x) + \delta,
    \\
    \label{eq:ApproximateSmoothness-RelationBetweenInnerProductAndSymmetricBregmanDistance}
    \InnerProduct{\bar{g}(x) - \bar{g}(y)}{x - y}
    \leq
    \beta_{f, \bar{f}, \bar{g}}(x, y) + \beta_{f, \bar{f}, \bar{g}}(y, x)
    \leq
    \InnerProduct{\bar{g}(x) - \bar{g}(y)}{x - y} + 2 \delta,
    \\
    \label{eq:ApproximateSmoothness-UpperBoundOnInnerProduct}
    \InnerProduct{\bar{g}(x) - \bar{g}(y)}{x - y}
    \leq
    L \Norm{x - y}^2 + 2 \delta,
    \\
    \label{eq:ApproximateSmoothness-LowerBoundOnBregmanDistance}
    \DualNorm{\bar{g}(x) - \bar{g}(y)}^2
    \leq
    2 L (\beta_{f, \bar{f}, \bar{g}}(x, y) + \delta),
    \\
    \label{eq:ApproximateSmoothness-UpperBoundOnMixedSquaredGradientNorm}
    \DualNorm{\Gradient f(x) - \bar{g}(y)}^2
    \leq
    2 L (\BregmanDistanceWithSubgradient{f}{\Gradient f(x)}(x, y) + \delta),
    \\
    \label{eq:ApproximateSmoothness-UpperBoundOnSquaredGradientDifference}
    \DualNorm{\bar{g}(x) - \bar{g}(y)}^2
    \leq
    L^2 \Norm{x - y}^2 + 4 L \delta,
    \\
    \label{eq:ApproximateSmoothness-UpperBoundOnGradientNormViaUsualBregmanDistance}
    \DualNorm{\bar{g}(x) - \bar{g}(y)}^2
    \leq
    4 L (\BregmanDistanceWithSubgradient{f}{\Gradient f(x)}(x, y) + 2 \delta),
    \\
    \label{eq:ApproximateSmoothness-UpperBoundOnUsualBregmanDistance}
    \BregmanDistanceWithSubgradient{f}{\Gradient f(x)}(x, y)
    \leq
    L \Norm{x - y}^2 + 2 \delta.
  \end{gather}
\end{theorem}

\begin{proof}
  Inequality~\eqref{eq:ApproximateSmoothness-RelationBetweenFunctionValues}
  follows immediately from
  \cref{eq:ApproximateSmoothness-UpperBoundOnBregmanDistance}
  by substituting $y = x$.

  To prove
  \cref{eq:ApproximateSmoothness-RelationBetweenInnerProductAndSymmetricBregmanDistance},
  we rewrite
  \[
    \beta_{f, \bar{f}, \bar{g}}(x, y) + \beta_{f, \bar{f}, \bar{g}}(y, x)
    =
    \InnerProduct{\bar{g}(x) - \bar{g}(y)}{x - y}
    +
    [f(x) - \bar{f}(x)] + [f(y) - \bar{f}(y)],
  \]
  and then apply \cref{eq:ApproximateSmoothness-RelationBetweenFunctionValues}.

  Using the first part of
  \cref{eq:ApproximateSmoothness-RelationBetweenInnerProductAndSymmetricBregmanDistance}
  and applying \cref{eq:ApproximateSmoothness-UpperBoundOnBregmanDistance} twice,
  we obtain \cref{eq:ApproximateSmoothness-UpperBoundOnInnerProduct}.

  To prove \cref{%
    eq:ApproximateSmoothness-LowerBoundOnBregmanDistance,%
    eq:ApproximateSmoothness-UpperBoundOnMixedSquaredGradientNorm%
  }, let us fix some $\bar{f}_1(x) \in \RealField$ and
  $\bar{g}_1(x) \in \RealField^d$ such that
  $
    \beta_{f, \bar{f}_1, \bar{g}_1}(z)
    \DefinedEqual
    f(z) - \bar{f}_1(x) - \InnerProduct{\bar{g}_1(x)}{z - x}
    \geq
    0
  $
  for any $z \in \RealField^d$.
  Note that we can choose either $(\bar{f}_1, \bar{g}_1) = (\bar{f}, \bar{g})$
  or $(\bar{f}_1, \bar{g}_1) = (f, \Gradient f)$.
  In view of \cref{eq:ApproximateSmoothness-UpperBoundOnBregmanDistance},
  for any $z \in \RealField^d$, we can write the following inequalities:
  \begin{align*}
    0
    \leq
    \beta_{f, \bar{f}_1, \bar{g}_1}(z)
    &\leq
    \bar{f}(y) - \bar{f}_1(x) - \InnerProduct{\bar{g}_1(x)}{y - x}
    +
    \InnerProduct{\bar{g}(y) - \bar{g}_1(x)}{z - y}
    +
    \frac{L}{2} \Norm{z - y}^2
    +
    \delta.
  \end{align*}
  Minimizing the right-hand side in $z \in \RealField^d$ and rearranging,
  we conclude that
  \[
    \frac{1}{2 L} \DualNorm{\bar{g}(y) - \bar{g}_1(x)}^2
    \leq
    \bar{f}(y) - \bar{f}_1(x) - \InnerProduct{\bar{g}_1(x)}{y - x} + \delta
    \leq
    \beta_{f, \bar{f}_1, \bar{g}_1}(x, y) + \delta,
  \]
  where the final inequality is due to
  \cref{eq:ApproximateSmoothness-RelationBetweenFunctionValues}.
  Substituting now either $(\bar{f}_1, \bar{g}_1) = (\bar{f}, \bar{g})$
  or $(\bar{f}_1, \bar{g}_1) = (f, \Gradient f)$, we obtain
  either \cref{eq:ApproximateSmoothness-LowerBoundOnBregmanDistance}
  or \cref{eq:ApproximateSmoothness-UpperBoundOnMixedSquaredGradientNorm},
  respectively.

  Inequality~\eqref{eq:ApproximateSmoothness-UpperBoundOnSquaredGradientDifference}
  follows immediately from \cref{%
    eq:ApproximateSmoothness-LowerBoundOnBregmanDistance,%
    eq:ApproximateSmoothness-UpperBoundOnBregmanDistance%
  }.

  Inequality~\eqref{eq:ApproximateSmoothness-UpperBoundOnGradientNormViaUsualBregmanDistance}
  follows from
  \cref{eq:ApproximateSmoothness-UpperBoundOnMixedSquaredGradientNorm}:
  \begin{align*}
    \DualNorm{\bar{g}(x) - \bar{g}(y)}^2
    &\leq
    2 \DualNorm{\Gradient f(x) - \bar{g}(y)}^2
    +
    2 \DualNorm{\bar{g}(x) - \Gradient f(x)}^2
    \\
    &\leq
    4 L (\BregmanDistanceWithSubgradient{f}{\Gradient f(x)}(x, y) + \delta)
    +
    4 L \delta
    =
    4 L (\BregmanDistanceWithSubgradient{f}{\Gradient f(x)}(x, y) + 2 \delta).
  \end{align*}

  To prove \cref{eq:ApproximateSmoothness-UpperBoundOnUsualBregmanDistance},
  we proceed as follows using
  first \cref{eq:ApproximateSmoothness-UpperBoundOnBregmanDistance},
  then \cref{eq:ApproximateSmoothness-RelationBetweenFunctionValues},
  and then \cref{eq:ApproximateSmoothness-UpperBoundOnMixedSquaredGradientNorm}:
  \begin{align*}
    \BregmanDistanceWithSubgradient{f}{\Gradient f(x)}(x, y)
    &\equiv
    f(y) - f(x) - \InnerProduct{\Gradient f(x)}{y - x}
    \\
    &\leq
    \bar{f}(x) - f(x)
    +
    \InnerProduct{\bar{g}(x) - \Gradient f(x)}{y - x}
    +
    \frac{L}{2} \Norm{y - x}^2 + \delta
    \\
    &\leq
    \InnerProduct{\bar{g}(x) - \Gradient f(x)}{y - x}
    +
    \frac{L}{2} \Norm{y - x}^2 + \delta
    \\
    &\leq
    \sqrt{2 L \delta} \Norm{y - x} + \frac{L}{2} \Norm{y - x}^2 + \delta
    \\
    &=
    \Bigl( \sqrt{\frac{L}{2}} \Norm{y - x} + \sqrt{\delta} \, \Bigr)^2
    \leq
    L \Norm{y - x}^2 + 2 \delta,
  \end{align*}
  where the final inequality is $(a + b)^2 \leq 2 a^2 + 2 b^2$, $a, b \geq 0$.
\end{proof}

\begin{remark}
  Some of the inequalities from
  \cref{th:InequalitiesForApproximatelySmoothFunctions}, namely, \cref{%
    eq:ApproximateSmoothness-RelationBetweenFunctionValues,%
    eq:ApproximateSmoothness-LowerBoundOnBregmanDistance,%
    eq:ApproximateSmoothness-UpperBoundOnSquaredGradientDifference%
  }, were established already in~\cite{devolder2013first}.
  We nevertheless prefer to present the corresponding proofs since they are
  rather simple, and we use the associated ideas for proving the other new
  inequalities.
\end{remark}

\begin{lemma}
  \label{th:StandardStochasticGradientOracleIsApproximatelySmooth}
  Let $\Map{f}{\RealField^d}{\RealField}$ be the function
  $f(x) \DefinedEqual \Expectation_{\xi}[f_{\xi}(x)]$,
  where each $\Map{f_{\xi}}{\RealField^d}{\RealField}$ is
  convex and $(\delta_{\xi}, L_{\xi})$-approximately smooth with components
  $(\bar{f}_{\xi}, \bar{g}_{\xi})$.
  Further, let $\hat{g}$ be the stochastic oracle defined by
  $g(x, \xi) \DefinedEqual \bar{g}_{\xi}(x)$,
  and let $\bar{f}(x) \DefinedEqual \Expectation_{\xi}[\bar{f}_{\xi}(x)]$,
  $\bar{g}(x) \DefinedEqual \Expectation_{\xi}[\bar{g}_{\xi}(x)]$.
  Then, $\hat{g}$ is an unbiased oracle for $\bar{g}$ and,
  for any $x, y \in \RealField^d$,
  $L_{\max} \DefinedEqual \sup_{\xi} L_{\xi}$ and
  $
    \bar{\delta}
    \DefinedEqual
    \frac{1}{L_{\max}} \Expectation_{\xi}[L_{\xi} \delta_{\xi}]
  $,
  it holds that
  \begin{equation}
    \label{eq:StandardStochasticGradientOracleIsApproximatelySmooth-Main}
    \Variance_{\hat{g}}(x, y)
    \leq
    2 L_{\max} [\beta_{f, \bar{f}, \bar{g}}(x, y) + \bar{\delta}].
  \end{equation}
  Furthermore, for any $x, y \in \RealField^d$ and
  any $\Gradient f(x) \in \Subdifferential f(x)$, it also holds that
  \begin{equation}
    \label{eq:StandardStochasticGradientOracleIsApproximatelySmooth-Extra}
    \Variance_{\hat{g}}(x, y)
    \leq
    4 L_{\max} [
      \BregmanDistanceWithSubgradient{f}{\Gradient f(x)}(x, y)
      +
      2 \bar{\delta}
    ].
  \end{equation}
\end{lemma}

\begin{proof}
  According to our definition of~$\bar{g}$, we have
  $\Expectation_{\xi}[\bar{g}_{\xi}(x)] = \bar{g}(x)$ for any $x$,
  so $\hat{g}$ is indeed an unbiased oracle for $\bar{g}$.
  Further, for any $x, y \in \RealField^d$, we can estimate
  \begin{align*}
    \Variance_{\hat{g}}(x, y)
    &\equiv
    \Expectation_{\xi}\bigl[
      \DualNorm{
        [\bar{g}_{\xi}(x) - \bar{g}_{\xi}(y)]
        -
        [\bar{g}(x) - \bar{g}(y)]
      }^2
    \bigr]
    \\
    &\leq
    \Expectation_{\xi} \bigl[
      \DualNorm{\bar{g}_{\xi}(x) - \bar{g}_{\xi}(y)}^2
    \bigr]
    \leq
    \Expectation_{\xi} \bigl[
      2 L_{\xi} \bigl(
        \beta_{f_{\xi}, \bar{f}_{\xi}, \bar{g}_{\xi}}(x, y) + \delta_{\xi}
      \bigr)
    \bigr]
    \\
    &\leq
    2 L_{\max} \bigl(
      \Expectation_{\xi}[\beta_{f_{\xi}, \bar{f}_{\xi}, \bar{g}_{\xi}}(x, y)]
      +
      \bar{\delta}
    \bigr)
    =
    2 L_{\max} [\beta_{f, \bar{f}, \bar{g}}(x, y) + \bar{\delta}],
  \end{align*}
  where $\bar{\delta}$ is as defined in the statement;
  the second inequality follows from
  \cref{th:InequalitiesForApproximatelySmoothFunctions}
  (inequality~\eqref{eq:ApproximateSmoothness-LowerBoundOnBregmanDistance}),
  and the final identity is due to the linearity of
  $\beta_{f, \bar{f}, \bar{g}}(x, y)$ in $(f, \bar{f}, \bar{g})$
  and the fact that, by our definitions,
  $\Expectation_{\xi}[f_{\xi}(x)] = f(x)$,
  $\Expectation_{\xi}[\bar{f}_{\xi}(x)] = \bar{f}(x)$,
  $\Expectation_{\xi}[\bar{g}_{\xi}(x)] = \bar{g}(x)$
  for any $x$.
  This proves
  \cref{eq:StandardStochasticGradientOracleIsApproximatelySmooth-Main}.

  The proof of
  \cref{eq:StandardStochasticGradientOracleIsApproximatelySmooth-Extra} is
  similar but now we apply
  \cref{eq:ApproximateSmoothness-UpperBoundOnGradientNormViaUsualBregmanDistance}
  instead of \cref{eq:ApproximateSmoothness-LowerBoundOnBregmanDistance}:
  \begin{align*}
    \Variance_{\hat{g}}(x, y)
    &\leq
    \Expectation_{\xi} \bigl[
      \DualNorm{\bar{g}_{\xi}(x) - \bar{g}_{\xi}(y)}^2
    \bigr]
    \leq
    \Expectation_{\xi} \bigl[
      4 L_{\xi} \bigl(
        \BregmanDistanceWithSubgradient{f_{\xi}}{\Gradient f_{\xi}(x)}(x, y)
        +
        2 \delta_{\xi}
      \bigr)
    \bigr]
    \\
    &\leq
    4 L_{\max} \bigl(
      \Expectation_{\xi}[
        \BregmanDistanceWithSubgradient{f_{\xi}}{\Gradient f_{\xi}(x)}(x, y)
      ]
      +
      2 \bar{\delta}
    \bigr)
    =
    4 L_{\max} [
      \BregmanDistanceWithSubgradient{f}{\Gradient f(x)}(x, y)
      +
      2 \bar{\delta}
    ],
  \end{align*}
  where we have used the fact that
  $\Subdifferential f(x) = \Expectation_{\xi}[\Subdifferential f_{\xi}(x)]$
  (see Proposition~2.2 in~\cite{bertsekas1973stochastic}),
  meaning that, for any $\Gradient f(x) \in \Subdifferential f(x)$, we can
  find a selection of $\Gradient f_{\xi}(x) \in \Subdifferential f_{\xi}(x)$
  such that $\Gradient f(x) = \Expectation_{\xi}[\Gradient f_{\xi}(x)]$.
\end{proof}
\subsection{Miscellaneous}

\begin{lemma}
  \label{th:OptimalityConditionForProximalStep}
  Let $\Map{\psi}{\RealField^d}{\RealFieldPlusInfty}$ be a proper closed convex
  function,
  $x \in \EffectiveDomain \psi$, $g \in \RealField^d$, $M \geq 0$,
  and let
  \[
    x_+ \DefinedEqual \ProximalMap_{\psi}(x, g, M).
  \]
  Then, for any $y \in \EffectiveDomain \psi$, we have
  \[
    \InnerProduct{g}{y - x_+} + \psi(y) + \frac{M}{2} \Norm{x - y}^2
    \geq
    \psi(x_+) + \frac{M}{2} \Norm{x - x_+}^2
    +
    \frac{M}{2} \Norm{x_+ - y}^2.
  \]
\end{lemma}

\begin{proof}
  Indeed, by definition, $x_+$ is the minimizer of the function
  $\Map{h}{\RealField^d}{\RealFieldPlusInfty}$ given by
  $
    h(y)
    \DefinedEqual
    \InnerProduct{g}{y} + \psi(y) + \frac{M}{2} \Norm{x - y}^2
  $,
  which is strongly convex with parameter~$M$ (or simply convex if $M = 0$).
  Hence, for any $y \in \EffectiveDomain \psi$ ($= \EffectiveDomain h$),
  we have $h(y) \geq h(x_+) + \frac{M}{2} \Norm{y - x_+}^2$,
  which is exactly the claimed inequality.
\end{proof}

\begin{lemma}
  \label{th:TelescopingDifferencesWithMin}
  Let $N \geq 1$ be an integer, $(M_k)_{k = 0}^N$ be a nondecreasing
  nonnegative sequence of reals, and let $\bar{M} \geq 0$.
  Then,
  \[
    \sum_{k = 0}^{N - 1}
    \PositivePart{\min\Set{M_{k + 1}, \bar{M}} - M_k}
    =
    \PositivePart{\min\Set{M_N, \bar{M}} - M_0}.
  \]
\end{lemma}

\begin{proof}
  It suffices to prove the identity only in the special case when $N = 2$,
  i.e., to show that $\gamma_0 + \gamma_1 = \Gamma$, where
  $\gamma_0 \DefinedEqual \PositivePart{\min\Set{M_1, \bar{M}} - M_0}$,
  $\gamma_1 \DefinedEqual \PositivePart{\min\Set{M_2, \bar{M}} - M_1}$,
  $\Gamma \DefinedEqual \PositivePart{\min\Set{M_2, \bar{M}} - M_0}$.
  The general case then easily follows by induction.

  To prove the identity, we use our assumption that $M_0 \leq M_1 \leq M_2$
  and consider three possible cases.
  If $M_1 \geq \bar{M}$, then
  $
    \gamma_0 + \gamma_1
    =
    \PositivePart{\bar{M} - M_0} + 0
    =
    \PositivePart{\bar{M} - M_0}
    =
    \Gamma
  $.
  If $M_1 < \bar{M} \leq M_2$, then
  $
    \gamma_0 + \gamma_1
    =
    (M_1 - M_0) + (\bar{M} - M_1)
    =
    \bar{M} - M_0
    =
    \Gamma
  $.
  Finally, if $M_2 < \bar{M}$, then
  $
    \gamma_0 + \gamma_1 = (M_1 - M_0) + (M_2 - M_1) = M_2 - M_0 = \Gamma
  $.
\end{proof}

\begin{lemma}
  \label{th:YoungsInequalityForVariance}
  Let $\hat{g}$ be a stochastic oracle in~$\RealField^d$.
  Then, for any $x, y, z \in \RealField^d$ and any $\tau > 0$, we have
  \begin{gather*}
    \Variance_{\hat{g}}(x)
    \leq
    (1 + \tau) \Variance_{\hat{g}}(y)
    +
    (1 + \tau^{-1}) \Variance_{\hat{g}}(x, y),
    \\
    \Variance_{\hat{g}}(x, y)
    \leq
    (1 + \tau) \Variance_{\hat{g}}(x, z)
    +
    (1 + \tau^{-1}) \Variance_{\hat{g}}(y, z).
  \end{gather*}
\end{lemma}

\begin{proof}
  Both inequalities are direct consequences of the standard inequality
  $
    \DualNorm{s_1 + s_2}^2
    \leq
    (1 + \tau) \DualNorm{s_1}^2 + (1 + \tau^{-1}) \DualNorm{s_2}^2
  $
  which is valid for any $s_1, s_2 \in \RealField^d$ and any $\tau > 0$.
  Indeed, let $g$ and $\xi$ be, respectively, the function and the random
  variable components of~$\hat{g}$, and let
  $\Delta(x, \xi) \DefinedEqual g(x, \xi) - \Expectation[g(x, \xi)]$
  for any $x \in \RealField^d$.
  Then, for any $x, y, z \in \RealField^d$ and $\tau > 0$, we can estimate
  \begin{align*}
    \Variance_{\hat{g}}(x)
    &\equiv
    \Expectation[\DualNorm{\Delta(x, \xi)}^2]
    =
    \Expectation[
      \DualNorm{\Delta(y, \xi) + [\Delta(x, \xi) - \Delta(y, \xi)]}^2
    ]
    \\
    &\leq
    (1 + \tau) \Expectation[\DualNorm{\Delta(y, \xi)}^2]
    +
    (1 + \tau^{-1}) \Expectation[\DualNorm{\Delta(x, \xi) - \Delta(y, \xi)}^2]
    \\
    &\equiv
    (1 + \tau) \Variance_{\hat{g}}(y)
    +
    (1 + \tau^{-1}) \Variance_{\hat{g}}(x, y).
  \end{align*}
  Similarly,
  \begin{align*}
    \Variance_{\hat{g}}(x, y)
    &\equiv
    \Expectation\bigl[ \DualNorm{\Delta(x, \xi) - \Delta(y, \xi)}^2 \bigr]
    =
    \Expectation\bigl[
      \DualNorm{
        [\Delta(x, \xi) - \Delta(z, \xi)] - [\Delta(y, \xi) - \Delta(z, \xi)]
      }^2
    \bigr]
    \\
    &\leq
    (1 + \tau)
    \Expectation\bigl[ \DualNorm{\Delta(x, \xi) - \Delta(z, \xi)}^2 \bigr]
    +
    (1 + \tau^{-1})
    \Expectation\bigl[ \DualNorm{\Delta(y, \xi) - \Delta(z, \xi)}^2 \bigr]
    \\
    &\equiv
    (1 + \tau) \Variance_{\hat{g}}(x, z)
    +
    (1 + \tau^{-1}) \Variance_{\hat{g}}(y, z).
    \qedhere
  \end{align*}
\end{proof}

\section{Omitted Proofs for \cref{sec:MainAlgorithmsAndStepsizeUpdateRules}}

\begin{lemma}[AdaGrad stepsize]
  \label{th:AdaGradRule}
  Let function~$f$ satisfy \cref{as:ApproximateSmoothness}.
  Consider the stepsize update rule
  $\hat{M}_+ = M_+(M, \Omega, x, \hat{x}_+, \hat{g}_x, \hat{g}_{x_+})$
  defined by
  \[
    \hat{M}_+
    \DefinedEqual
    \sqrt{M^2 + \frac{1}{\Omega} \DualNorm{\hat{g}_{x_+} - \hat{g}_x}^2}.
  \]
  Then, this stepsize update rules satisfies
  \cref{eq:RequirementOnStepsizeUpdateRule}
  with $c_1 = \frac{5}{2}$, $c_2 = 4$, $c_3 = 6$, $c_4 = 2$.
\end{lemma}

\begin{proof}
  Let
  $
    \hat{\Delta}(\bar{M})
    \DefinedEqual
    \beta_{f, \bar{f}, \bar{g}}(x, \hat{x}_+)
    +
    \InnerProduct{\bar{g}(x) - \hat{g}_x}{\hat{x}_+ - x}
    -
    \frac{\bar{M}}{2} \Norm{\hat{x}_+ - x}^2
  $.
  From our \cref{as:ApproximateSmoothness} and
  \cref{th:InequalitiesForApproximatelySmoothFunctions}
  (inequality~\eqref{eq:ApproximateSmoothness-RelationBetweenInnerProductAndSymmetricBregmanDistance}),
  it follows that
  $
    \beta_{f, \bar{f}, \bar{g}}(x, \hat{x}_+)
    +
    \beta_{f, \bar{f}, \bar{g}}(\hat{x}_+, x)
    \leq
    \InnerProduct{\bar{g}(\hat{x}_+) - \bar{g}(x)}{\hat{x}_+ - x}
    +
    2 \delta_f
  $.
  Hence,
  \[
    \Expectation[
      \hat{\Delta}(\hat{M}_+) + \beta_{f, \bar{f}, \bar{g}}(\hat{x}_+, x)
    ]
    \leq
    \Expectation\Bigl[
      \InnerProduct{\bar{g}(\hat{x}_+) - \hat{g}_x}{\hat{x}_+ - x}
      -
      \frac{\hat{M}_+}{2} \Norm{\hat{x}_+ - x}^2
    \Bigr]
    +
    2 \delta{_f}
    =
    \Expectation[\hat{\Delta}_1(\hat{M}_+)] + 2 \delta_f,
  \]
  where
  $
    \hat{\Delta}_1(\hat{M}_+)
    \DefinedEqual
    \InnerProduct{\hat{g}_{x_+} - \hat{g}_x}{\hat{x}_+ - x}
    -
    \frac{\hat{M}_+}{2} \Norm{\hat{x}_+ - x}^2
  $.
  Hence,
  \[
    \Gamma
    \DefinedEqual
    \Expectation[
      \hat{\Delta}(\hat{M}_+) + (\hat{M}_+ - M) \Omega
      +
      \beta_{f, \bar{f}, \bar{g}}(\hat{x}_+, x)
    ]
    \leq
    \Expectation[\hat{\Delta}_1(\hat{M}_+) + (\hat{M}_+ - M) \Omega]
    +
    2 \delta_f.
  \]

  From the definition of~$\hat{M}_+$, it follows that
  $
    \DualNorm{\hat{g}_{x_+} - \hat{g}_x}^2
    =
    (\hat{M}_+^2 - M^2) \Omega
    =
    (\hat{M}_+ + M) (\hat{M}_+ - M) \Omega.
  $
  Since $\hat{M}_+ \geq M$, this means that
  \[
    \frac{1}{2 \hat{M}_+} \DualNorm{\hat{g}_{x_+} - \hat{g}_x}^2
    \leq
    (\hat{M}_+ - M) \Omega
    \leq
    \frac{1}{\hat{M}_+} \DualNorm{\hat{g}_{x_+} - \hat{g}_x}^2
  \]
  Let us now upper bound
  $
    \hat{\Gamma}
    \DefinedEqual
    \hat{\Delta}_1(\hat{M}_+) + (\hat{M}_+ - M) \Omega
  $.
  For this, let us fix an arbitrary constant $\bar{M} \geq 0$ and consider
  two cases.
  If $\hat{M}_+ \geq \bar{M}$, we can bound
  \[
    \hat{\Gamma}
    \leq
    \hat{\Delta}_1(\hat{M}_+)
    +
    \frac{1}{\hat{M}_+} \DualNorm{\hat{g}_{x_+} - \hat{g}_x}^2
    \leq
    \hat{\Delta}_1(\bar{M})
    +
    \frac{1}{\bar{M}} \DualNorm{\hat{g}_{x_+} - \hat{g}_x}^2
    \EqualDefines
    \hat{\Gamma}(\bar{M}).
  \]
  If $\hat{M}_+ \leq \bar{M}$, we can bound
  \[
    \hat{\Gamma}
    \leq
    \frac{1}{2 \hat{M}_+} \DualNorm{\hat{g}_{x_+} - \hat{g}_x}^2
    +
    (\hat{M}_+ - M) \Omega
    \leq
    2 (\hat{M}_+ - M) \Omega
    =
    2 \PositivePart{\min\Set{\hat{M}_+, \bar{M}} - M} \Omega.
  \]
  Combining the two cases, we get
  $
    \hat{\Gamma}
    \leq
    \PositivePart{\hat{\Gamma}(\bar{M})}
    +
    2 \PositivePart{\min\Set{\hat{M}_+, \bar{M}} - M} \Omega
  $.
  Thus,
  \[
    \Gamma
    \leq
    \Expectation[\hat{\Gamma}] + 2 \delta_f
    \leq
    \Expectation\bigl\{ \PositivePart{\hat{\Gamma}(\bar{M})} \bigr\}
    +
    2 \Expectation\bigl\{
      \PositivePart{\min\Set{\hat{M}_+, \bar{M}} - M} \Omega
    \bigr\}
    +
    2 \delta_f.
  \]

  Let us now estimate the first term.
  Denote $\hat{S} \DefinedEqual \hat{g}_x - \bar{g}(x)$ and
  $\hat{S}_+ \DefinedEqual \hat{g}_{x_+} - \bar{g}(\hat{x}_+)$.
  Then,
  \begin{align*}
    \hat{\Gamma}(\bar{M})
    &\equiv
    \InnerProduct{\hat{g}_{x_+} - \hat{g}_x}{\hat{x}_+ - x}
    -
    \frac{\bar{M}}{2} \Norm{\hat{x}_+ - x}^2
    +
    \frac{1}{\bar{M}} \DualNorm{\hat{g}_{x_+} - \hat{g}_x}^2
    \\
    &\leq
    \InnerProduct{\bar{g}(\hat{x}_+) - \bar{g}(x)}{\hat{x}_+ - x}
    +
    \frac{2}{\bar{M}} \DualNorm{\bar{g}(\hat{x}_+) - \bar{g}(x)}^2
    \\
    &\qquad
    +
    \InnerProduct{\hat{S}_+ - \hat{S}}{\hat{x}_+ - x}
    +
    \frac{2}{\bar{M}} \DualNorm{\hat{S}_+ - \hat{S}}^2
    -
    \frac{\bar{M}}{2} \Norm{\hat{x}_+ - x}^2
  \end{align*}
  Using now our \cref{as:ApproximateSmoothness} and
  \cref{th:InequalitiesForApproximatelySmoothFunctions}
  (inequalities~\eqref{eq:ApproximateSmoothness-UpperBoundOnInnerProduct}
    and~\eqref{eq:ApproximateSmoothness-UpperBoundOnSquaredGradientDifference}),
  we can continue as follows:
  \begin{align*}
    \hat{\Gamma}(\bar{M})
    &\leq
    L_f \Norm{\hat{x}_+ - x}^2 + 2 \delta_f
    +
    \frac{2}{\bar{M}} (L_f^2 \Norm{\hat{x}_+ - x}^2 + 4 L_f \delta_f)
    \\
    &\qquad
    +
    \InnerProduct{\hat{S}_+ - \hat{S}}{\hat{x}_+ - x}
    +
    \frac{2}{\bar{M}} \DualNorm{\hat{S}_+ - \hat{S}}^2
    -
    \frac{\bar{M}}{2} \Norm{\hat{x}_+ - x}^2
    \\
    &\leq
    \InnerProduct{\hat{S}_+ - \hat{S}}{\hat{x}_+ - x}
    +
    \frac{2}{\bar{M}} \DualNorm{\hat{S}_+ - \hat{S}}^2
    -
    \frac{\bar{M} - 2 L_f (1 + \frac{2 L_f}{\bar{M}})}{2}
    \Norm{\hat{x}_+ - x}^2
    +
    2 \Bigl( 1 + \frac{4 L_f}{\bar{M}} \Bigr) \delta_f
    \\
    &\leq
    \Bigl(
      \frac{2}{\bar{M}}
      +
      \frac{1}{2 [\bar{M} - 2 L_f (1 + \frac{2 L_f}{\bar{M}})]}
    \Bigr) \DualNorm{\hat{S}_+ - \hat{S}}^2
    +
    2 \Bigl( 1 + \frac{4 L_f}{\bar{M}} \Bigr) \delta_f.
  \end{align*}
  Consequently,
  \[
    \Expectation\bigl\{ \PositivePart{\hat{\Gamma}(\bar{M})} \bigr\}
    \leq
    \Bigl(
      \frac{2}{\bar{M}}
      +
      \frac{1}{2 [\bar{M} - 2 L_f (1 + \frac{2 L_f}{\bar{M}})]}
    \Bigr) \Expectation[\DualNorm{\hat{S}_+ - \hat{S}}^2]
    +
    2 \Bigl( 1 + \frac{4 L_f}{\bar{M}} \Bigr) \delta_f.
  \]
  In particular, for $\bar{M} > 4 L_f$, we can estimate
  $
    \frac{2}{\bar{M}}
    +
    \frac{1}{2 [\bar{M} - 2 L_f (1 + \frac{2 L_f}{\bar{M}})]}
    \leq
    \frac{2}{\bar{M}}
    +
    \frac{1}{2 (\bar{M} - 4 L_f)}
    \leq
    \frac{5}{2 (\bar{M} - 4 L_f)}.
  $
  Therefore, for any $\bar{M} > 4 L_f$,
  \[
    \Expectation\bigl\{ \PositivePart{\hat{\Gamma}(\bar{M})} \bigr\}
    \leq
    \frac{5}{2 (\bar{M} - 4 L_f)} \Expectation[\DualNorm{\hat{S}_+ - \hat{S}}^2]
    +
    4 \delta_f
    =
    \frac{5}{2 (\bar{M} - 4 L_f)}
    \Expectation[\Variance_{\hat{g}}(\hat{x}_+) + \Variance_{\hat{g}}(x)]
    +
    4 \delta_f,
  \]
  where the final identity follows from the fact that
  $
    \Expectation[\DualNorm{\hat{S}_+ - \hat{S}}^2]
    =
    \Expectation[\DualNorm{\hat{S}_+}^2]
    +
    \Expectation[\DualNorm{\hat{S}}^2]
    =
    \Expectation[\Variance_{\hat{g}}(\hat{x}_+)]
    +
    \Variance_{\hat{g}}(x)
  $
  (because $\hat{S}_+$, conditioned on the randomness coming from~$\hat{g}_x$,
    has zero mean).

  Combining everything together, we get
  \[
    \Gamma
    \leq
    \frac{5}{2 (\bar{M} - 4 L_f)}
    \Expectation[\Variance_{\hat{g}}(\hat{x}_+) + \Variance_{\hat{g}}(x)]
    +
    6 \delta_f
    +
    2 \Expectation\bigl\{
      \PositivePart{\min\Set{\hat{M}_+, \bar{M}} - M} \Omega
    \bigr\}.
  \]
  This is exactly \cref{eq:RequirementOnStepsizeUpdateRule} with
  $c_1 = \frac{5}{2}$, $c_2 = 4$, $c_3 = 6$, $c_4 = 2$.
\end{proof}

\begin{lemma}
  \label{th:ModifiedAdaGradRule}
  Let function~$f$ satisfy \cref{as:ApproximateSmoothness}.
  Consider the stepsize update rule
  $\hat{M}_+ = M_+(M, \Omega, x, \hat{x}_+, \hat{g}_x, \hat{g}_{x_+})$
  defined as the solution of the following equation:
  \[
    (\hat{M}_+ - M) \Omega = \PositivePart{\hat{\Delta}_1(\hat{M}_+)},
    \qquad
    \hat{\Delta}_1(\hat{M}_+)
    \DefinedEqual
    \InnerProduct{\hat{g}_{x_+} - \hat{g}_x}{\hat{x}_+ - x}
    -
    \frac{\hat{M}_+}{2} \Norm{\hat{x}_+ - x}^2.
  \]
  Then, this stepsize update rules satisfies
  \cref{eq:RequirementOnStepsizeUpdateRule}
  with $c_1 = 1$, $c_2 = 2$, $c_3 = 6$, $c_4 = 2$.
\end{lemma}

\begin{proof}
  Let us define
  $
    \hat{\Delta}(\bar{M})
    \DefinedEqual
    \beta_{f, \bar{f}, \bar{g}}(x, \hat{x}_+)
    +
    \InnerProduct{\bar{g}(x) - \hat{g}_x}{\hat{x}_+ - x}
    -
    \frac{\bar{M}}{2} \Norm{\hat{x}_+ - x}^2
  $.
  Starting as in the proof of \cref{th:AdaGradRule}, we see that
  \[
    \Gamma
    \DefinedEqual
    \Expectation[
      \hat{\Delta}(\hat{M}_+) + (\hat{M}_+ - M) \Omega
      +
      \beta_{f, \bar{f}, \bar{g}}(\hat{x}_+, x)
    ]
    \leq
    \Expectation[\hat{\Delta}_1(\hat{M}_+) + (\hat{M}_+ - M) \Omega]
    +
    2 \delta_f,
  \]
  with the same $\hat{\Delta}_1(\cdot)$ as defined in the statement.

  Let us now upper bound
  $
    \hat{\Gamma}
    \DefinedEqual
    \hat{\Delta}_1(\hat{M}_+) + (\hat{M}_+ - M) \Omega
  $.
  For this, let us fix an arbitrary constant $\bar{M} \geq 0$ and consider
  two cases.
  If $\hat{M}_+ \geq \bar{M}$, we can bound, using the monotonicity of
  $\hat{\Delta}_1(\cdot)$,
  \[
    \hat{\Gamma}
    =
    \hat{\Delta}_1(\hat{M}_+)
    +
    \PositivePart{\hat{\Delta}_1(\hat{M}_+)}
    \leq
    \hat{\Delta}_1(\bar{M})
    +
    \PositivePart{\hat{\Delta}_1(\bar{M})}
    \leq
    2 \PositivePart{\hat{\Delta}_1(\bar{M})}.
  \]
  If $\hat{M}_+ \leq \bar{M}$, we can bound
  \[
    \hat{\Gamma}
    \leq
    \PositivePart{\hat{\Delta}_1(\hat{M}_+)} + (\hat{M}_+ - M) \Omega
    =
    2 (\hat{M}_+ - M) \Omega
    =
    2 \PositivePart{\min\Set{\hat{M}_+, \bar{M}} - M} \Omega.
  \]
  Combining the two cases, we get
  $
    \hat{\Gamma}
    \leq
    2 \PositivePart{\hat{\Delta}_1(\bar{M})}
    +
    2 \PositivePart{\min\Set{\hat{M}_+, \bar{M}} - M} \Omega,
  $
  and hence
  \[
    \Gamma
    \leq
    \Expectation[\hat{\Gamma}] + 2 \delta_f
    \leq
    2 \Expectation\bigl\{ \PositivePart{\hat{\Delta}_1(\bar{M})} \bigr\}
    +
    2 \Expectation\bigl\{
      \PositivePart{\min\Set{\hat{M}_+, \bar{M}} - M} \Omega
    \bigr\}
    +
    2 \delta_f.
  \]

  Let us now estimate the first term.
  According to our \cref{as:ApproximateSmoothness} and
  \cref{th:InequalitiesForApproximatelySmoothFunctions}
  (inequality~\eqref{eq:ApproximateSmoothness-UpperBoundOnInnerProduct}),
  we have
  $
    \InnerProduct{\bar{g}(\hat{x}_+) - \bar{g}(x)}{\hat{x}_+ - x}
    \leq
    L_f \Norm{\hat{x}_+ - x}^2 + 2 \delta_f
  $.
  Hence, denoting $\hat{S} \DefinedEqual \hat{g}_x - \bar{g}(x)$ and
  $\hat{S}_+ \DefinedEqual \hat{g}_{x_+} - \bar{g}(\hat{x}_+)$,
  we can estimate, for any $\bar{M} > 2 L_f$,
  \begin{align*}
    \hat{\Delta}_1(\bar{M})
    &=
    \InnerProduct{\bar{g}(\hat{x}_+) - \bar{g}(x)}{\hat{x}_+ - x}
    +
    \InnerProduct{\hat{S}_+ - \hat{S}}{\hat{x}_+ - x}
    -
    \frac{\bar{M}}{2} \Norm{\hat{x}_+ - x}^2
    \\
    &\leq
    \InnerProduct{\hat{S}_+ - \hat{S}}{\hat{x}_+ - x}
    -
    \frac{\bar{M} - 2 L_f}{2} \Norm{\hat{x}_+ - x}^2
    +
    2 \delta_f
    \leq
    \frac{1}{2 (\bar{M} - 2 L_f)} \DualNorm{\hat{S}_+ - \hat{S}}^2
    +
    2 \delta_f.
  \end{align*}
  Hence,
  \[
    \Expectation\bigl\{ \PositivePart{\hat{\Delta}_1(\bar{M})} \bigr\}
    \leq
    \frac{1}{2 (\bar{M} - 2 L_f)}
    \Expectation[\DualNorm{\hat{S}_+ - \hat{S}}^2]
    +
    2 \delta_f
    =
    \frac{1}{2 (\bar{M} - 2 L_f)} \Expectation[
      \Variance_{\hat{g}}(\hat{x}_+)
      +
      \Variance_{\hat{g}}(x)
    ]
    +
    2 \delta_f,
  \]
  where the final identity follows from the fact that
  $
    \Expectation[\DualNorm{\hat{S}_+ - \hat{S}}^2]
    =
    \Expectation[\DualNorm{\hat{S}_+}^2]
    +
    \Expectation[\DualNorm{\hat{S}}^2]
    =
    \Expectation[\Variance_{\hat{g}}(\hat{x}_+)]
    +
    \Variance_{\hat{g}}(x)
  $
  (because $\hat{S}_+$, conditioned on the randomness coming from~$\hat{g}_x$,
    has zero mean).

  Thus,
  \[
    \Gamma
    \leq
    \frac{1}{\bar{M} - 2 L_f}
    \Expectation[\Variance_{\hat{g}}(\hat{x}_+) + \Variance_{\hat{g}}(x)]
    +
    6 \delta_f
    +
    2 \Expectation\bigl\{
      \PositivePart{\min\Set{\hat{M}_+, \bar{M}} - M} \Omega
    \bigr\},
  \]
  which is exactly \cref{eq:RequirementOnStepsizeUpdateRule} with
  $c_1 = 1$, $c_2 = 2$, $c_3 = 6$, $c_4 = 2$.
\end{proof}

\section{Omitted Proofs for \cref{sec:UniformlyBoundedVariance}}
\label{sec:UniformlyBoundedVariance-Proofs}

\subsection{Universal SGD}
\label{sec:UniversalSgd-VarianceAtSolution-Proofs}

\thUniversalSgdVarianceAtSolution*

\begin{proof}
  \label{th:UniversalSgd-VarianceAtSolution:Proof}
  Let $x_0, \ldots, x_N$ be the points generated inside the method
  and let $F_N \DefinedEqual \Expectation[F(\bar{x}_N)] - F^*$.
  Using \cref{th:YoungsInequalityForVariance,as:ApproximatelySmoothVariance},
  we can estimate, for any $0 \leq k \leq N - 1$,
  \begin{align*}
    \Variance_{\hat{g}}(x_{k + 1}) + \Variance_{\hat{g}}(x_k)
    &\leq
    3 \Variance_{\hat{g}}(x_k) + 2 \Variance_{\hat{g}}(x_{k + 1}, x_k)
    \\
    &\leq
    6 \sigma_*^2
    +
    6 \Variance_{\hat{g}}(x_k, x^*)
    +
    2 \Variance_{\hat{g}}(x_{k + 1}, x_k)
    \\
    &\leq
    6 \sigma_*^2
    +
    12 L_{\hat{g}}[\beta_{f, \bar{f}, \bar{g}}(x_k, x^*) + \delta_{\hat{g}}]
    +
    4 L_{\hat{g}}[\beta_{f, \bar{f}, \bar{g}}(x_{k + 1}, x_k) + \delta_{\hat{g}}]
    \\
    &=
    6 \sigma_*^2
    +
    4 L_{\hat{g}} [
      3 \beta_{f, \bar{f}, \bar{g}}(x_k, x^*)
      +
      \beta_{f, \bar{f}, \bar{g}}(x_{k + 1}, x_k)
      +
      4 \delta_{\hat{g}}
    ].
  \end{align*}
  Substituting this bound into the general guarantee given by
  \cref{th:UniversalSgd-General}
  (and taking into account the fact that $M_0 = 0$), we obtain
  \begin{multline*}
    N F_N
    +
    \sum_{k = 0}^{N - 1} \Expectation[
      \beta_{f, \bar{f}, \bar{g}}(x_{k + 1}, x_k)
      +
      \beta_{f, \bar{f}, \bar{g}}(x_k, x^*)
    ]
    \\
    \leq
    c_4 \bar{M} D^2
    +
    \frac{6 c_1 \sigma_*^2 N}{\bar{M} - c_2 L_f}
    +
    \alpha \sum_{k = 0}^{N - 1} \Expectation[
      \beta_{f, \bar{f}, \bar{g}}(x_{k + 1}, x_k)
      +
      3 \beta_{f, \bar{f}, \bar{g}}(x_k, x^*)
    ]
    +
    N (c_3 \delta_f + 4 \alpha \delta_{\hat{g}}),
  \end{multline*}
  where $\bar{M} > c_2 L_f$ is an arbitrary constant
  and $\alpha \DefinedEqual \frac{4 c_1 L_{\hat{g}}}{\bar{M} - c_2 L_f}$.
  Requiring now that $3 \alpha \leq 1$ or, equivalently, that
  $\bar{M} \geq c_2 L_f + 12 c_1 L_{\hat{g}} \EqualDefines \bar{M}_{\min}$,
  we can cancel the nonnegative $\beta_{f, \bar{f}, \bar{g}}(\cdot, \cdot)$
  terms on both sides and obtain
  \[
    F_N
    \leq
    \frac{c_4 \bar{M} D^2}{N}
    +
    \frac{6 c_1 \sigma_*^2}{\bar{M} - c_2 L_f}
    +
    \delta,
  \]
  where $\delta \DefinedEqual c_3 \delta_f + \frac{4}{3} \delta_{\hat{g}}$.
  The optimal coefficient $\bar{M}_*$ minimizing the right-hand side is
  $\bar{M}_* = c_2 L_f + \frac{\sigma_*}{D} \sqrt{\frac{6 c_1 N}{c_4}}$.
  However, we still need to respect the constraint
  $\bar{M} \geq \bar{M}_{\min}$.
  Choosing
  $
    \bar{M}
    =
    c_2 L_f
    +
    12 c_1 L_{\hat{g}}
    +
    \frac{\sigma_*}{D} \sqrt{\frac{6 c_1 N}{c_4}}
  $,
  we conclude that
  \begin{align*}
    F_N
    &\leq
    \frac{c_4 D^2}{N} \Bigl(
      c_2 L_f
      +
      12 c_1 L_{\hat{g}}
      +
      \frac{\sigma_*}{D} \sqrt{\frac{6 c_1 N}{c_4}}
    \Bigr)
    +
    \frac{6 c_1 \sigma_*^2}{\frac{\sigma_*}{D} \sqrt{\frac{6 c_1 N}{c_4}}}
    +
    \delta
    \\
    &=
    \frac{c_4 (c_2 L_f + 12 c_1 L_{\hat{g}}) D^2}{N}
    +
    2 \sigma_* D \sqrt{\frac{6 c_1 c_4}{N}}
    +
    \delta.
    \qedhere
  \end{align*}
\end{proof}

\subsection{Universal Fast SGD}
\label{sec:UniversalFastSgd-VarianceAtSolution-Proofs}

\thUniversalFastSgdVarianceAtSolution*

\begin{proof}
  \label{th:UniversalFastSgd-VarianceAtSolution:Proof}
  Let $k \geq 1$ be arbitrary and
  $F_k \DefinedEqual \Expectation[F(x_k)] - F^*$.
  Using \cref{th:YoungsInequalityForVariance,as:ApproximatelySmoothVariance},
  we can estimate, for each $i$,
  \begin{align*}
    \Variance_{\hat{g}}(x_{i + 1}) + \Variance_{\hat{g}}(y_i)
    &\leq
    3 \Variance_{\hat{g}}(y_i) + 2 \Variance_{\hat{g}}(x_{i + 1}, y_i)
    \\
    &\leq
    6 \sigma_*^2
    +
    6 \Variance_{\hat{g}}(y_i, x^*)
    +
    2 \Variance_{\hat{g}}(x_{i + 1}, y_i)
    \\
    &\leq
    6 \sigma_*^2
    +
    12 L_{\hat{g}} [\beta_{f, \bar{f}, \bar{g}}(y_i, x^*) + \delta_{\hat{g}}]
    +
    4 L_{\hat{g}} [\beta_{f, \bar{f}, \bar{g}}(x_{i + 1}, y_i) + \delta_{\hat{g}}]
    \\
    &=
    6 \sigma_*^2
    +
    4 L_{\hat{g}} [
      3 \beta_{f, \bar{f}, \bar{g}}(y_i, x^*)
      +
      \beta_{f, \bar{f}, \bar{g}}(x_{i + 1}, y_i)
      +
      4 \delta_{\hat{g}}
    ].
  \end{align*}
  Substituting this bound into the guarantee given by
  \cref{th:UniversalFastSgd-GeneralGuarantee}, we obtain
  \begin{multline*}
    A_k F_k
    +
    \sum_{i = 0}^{k - 1} \Expectation[
      A_{i + 1} \beta_{f, \bar{f}, \bar{g}}(x_{i + 1}, y_i)
      +
      a_{i + 1} \beta_{f, \bar{f}, \bar{g}}(y_i, x^*)
    ]
    \\
    \leq
    c_4 \bar{M} D^2
    +
    \sum_{i = 0}^{k - 1} \alpha_{i + 1} \Expectation[
      \beta_{f, \bar{f}, \bar{g}}(x_{i + 1}, y_i)
      +
      3 \beta_{f, \bar{f}, \bar{g}}(y_i, x^*)
      +
      4 \delta_{\hat{g}}
    ]
    +
    \frac{6 c_1 \sigma_*^2}{\bar{M} - c_2 L_f} \sum_{i = 1}^k a_i^2
    +
    c_3 \delta_f \sum_{i = 1}^k A_i,
  \end{multline*}
  where
  $
    \alpha_{i + 1}
    \DefinedEqual
    \frac{4 c_1 L_{\hat{g}} a_{i + 1}^2}{\bar{M} - c_2 L_f}
  $,
  $a_i = \frac{1}{2} i$, $A_k = \frac{1}{4} k (k + 1)$,
  $\sum_{i = 1}^k a_i^2 = \frac{1}{24} k (k + 1) (2 k + 1)$,
  $\sum_{i = 1}^k A_i = \frac{1}{12} k (k + 1) (k + 2)$.
  Requiring now that $3 \alpha_{i + 1} \leq a_{i + 1}$
  for all $i = 0, \ldots, k - 1$ or, equivalently, that
  $
    \bar{M}
    \geq
    c_2 L_f + 12 c_1 L_{\hat{g}} a_k
    \equiv
    c_2 L_f + 6 c_1 L_{\hat{g}} k
  $,
  we can cancel the nonnegative $\beta_{f, \bar{f}, \bar{g}}(\cdot, \cdot)$
  terms on both sides and obtain
  \begin{align*}
    F_k
    &\leq
    \frac{1}{A_k} \Bigl(
      c_4 \bar{M} D^2
      +
      \frac{6 c_1 \sigma_*^2}{\bar{M} - c_2 L_f}
      \sum_{i = 1}^k a_i^2
      +
      c_3 \delta_f \sum_{i = 1}^k A_i
      +
      \frac{4}{3} A_k \delta_{\hat{g}}
    \Bigr)
    \\
    &=
    \frac{4}{k (k + 1)} \Bigl(
      c_4 \bar{M} D^2
      +
      \frac{c_1 \sigma_*^2 k (k + 1) (2 k + 1)}{4 (\bar{M} - c_2 L_f)}
      +
      \frac{c_3}{12} \delta_f k (k + 1) (k + 2)
    \Bigr)
    +
    \frac{4}{3} \delta_{\hat{g}}
    \\
    &=
    \frac{4 c_4 \bar{M} D^2}{k (k + 1)}
    +
    \frac{c_1 \sigma_*^2 (2 k + 1)}{\bar{M} - c_2 L_f}
    +
    \delta_k,
  \end{align*}
  where
  $
    \delta_k
    \DefinedEqual
    \frac{c_3}{3} (k + 2) \delta_f + \frac{4}{3} \delta_{\hat{g}}
  $.

  The minimizer of the right-hand side is
  $
    \bar{M}_*
    =
    c_2 L_f
    +
    \frac{\sigma_*}{2 D} \sqrt{\frac{c_1}{c_4} k (k + 1) (2 k + 1)}
  $.
  However, recall that we also need to satisfy the constraint
  $\bar{M} \geq c_2 L_f + 6 c_1 L_{\hat{g}} k$.
  Choosing
  $
    \bar{M}
    =
    c_2 L_f
    +
    6 c_1 L_{\hat{g}} k
    +
    \frac{\sigma_*}{2 D} \sqrt{\frac{c_1}{c_4} k (k + 1) (2 k + 1)}
  $,
  we obtain
  \begin{align*}
    F_k
    &\leq
    \begin{multlined}[t]
      \frac{4 c_4 D^2}{k (k + 1)} \Bigl(
        c_2 L_f
        +
        6 c_1 L_{\hat{g}} k
        +
        \frac{\sigma_*}{2 D} \sqrt{\frac{c_1}{c_4} k (k + 1) (2 k + 1)}
      \Bigr)
      +
      \frac{
        c_1 \sigma_*^2 (2 k + 1)
      }{
        \frac{\sigma_*}{2 D} \sqrt{\frac{c_1}{c_4} k (k + 1) (2 k + 1)}
      }
      +
      \delta_k
    \end{multlined}
    \\
    &=
    \frac{4 c_2 c_4 L_f D^2}{k (k + 1)}
    +
    \frac{24 c_1 c_4 L_{\hat{g}} D^2}{k + 1}
    +
    4 \sigma_* D \sqrt{\frac{c_1 c_4 (2 k + 1)}{k (k + 1)}}
    +
    \delta_k
    \\
    &\leq
    \frac{4 c_2 c_4 L_f D^2}{k (k + 1)}
    +
    \frac{24 c_1 c_4 L_{\hat{g}} D^2}{k + 1}
    +
    4 \sigma_* D \sqrt{\frac{2 c_1 c_4}{k}}
    +
    \delta_k.
    \qedhere
  \end{align*}
\end{proof}

\section{Omitted Proofs for \cref{sec:ImplicitVarianceReduction}}

\section{Omitted Proofs for \cref{sec:ExplicitVarianceReductionWithSvrg}}
\label{sec:ExplicitVarianceReductionWithSvrg-Proofs}

\begin{lemma}[Basic property of SVRG oracle]
  \label{th:SvrgOracle}
  Let $\hat{g}$ be a stochastic oracle in~$\RealField^d$, and let
  $\hat{G} = \SvrgOracle_{\hat{g}}(\tilde{x})$
  for some $\tilde{x} \in \RealField^d$.
  Then, for any $x \in \RealField^d$, the mean value of $\hat{G}$ at~$x$
  is the same as that of $\hat{g}$ at~$x$, while
  $\Variance_{\hat{G}}(x) = \Variance_{\hat{g}}(x, \tilde{x})$.
\end{lemma}

\begin{proof}
  Let $g$ and $\xi$ be, respectively, the function and the random variable
  components of~$\hat{g}$, and let
  $g(x) \DefinedEqual \Expectation_{\xi}[g(x, \xi)]$,
  $g(\tilde{x}) \DefinedEqual \Expectation_{\xi}[g(\tilde{x}, \xi)]$.
  Then, by definition, $\hat{G}$ is the oracle with the same random variable
  component~$\xi$ and the function component~$G$ defined by
  $G(x, \xi) = g(x, \xi) - g(\tilde{x}, \xi) + g(\tilde{x})$.
  Consequently, $\Expectation_{\xi}[G(x, \xi)] = g(x)$, and
  \begin{align*}
    \Variance_{\hat{G}}(x)
    &=
    \Expectation_{\xi}\bigl[ \DualNorm{G(x, \xi) - g(x)]}^2 \bigr]
    \\
    &=
    \Expectation_{\xi}\bigl[
      \DualNorm{[g(x, \xi) - g(\tilde{x}, \xi)] - [g(x) - g(\tilde{x})]}^2
    \bigr]
    =
    \Variance_{\hat{g}}(x, \tilde{x}).
    \qedhere
  \end{align*}
\end{proof}

\subsection{Universal SVRG}

\begin{lemma}[Universal SVRG Epoch]
  \label{th:UniversalSvrgEpoch}
  Consider problem~\eqref{eq:Problem} under
  \cref{%
    as:BoundedFeasibleSet,%
    as:ApproximateSmoothness,%
    as:ApproximatelySmoothVariance,%
    as:ApproximatelySmoothVariance-Extra%
  }.
  Let $x, \tilde{x} \in \EffectiveDomain \psi$ be points,
  $M \geq 0$ be a coefficient, $N \geq 1$ be an integer,
  $\hat{G} = \SvrgOracle_{\hat{g}}(\tilde{x})$, and let
  \[
    (\tilde{x}_+, x_+, M_+)
    \EqualRandom
    \UniversalSgd_{\hat{G}, \psi}(x, M, N; D),
  \]
  as defined by \cref{alg:UniversalSgd}.
  Then, for any $\bar{M} \geq c_2 L_f + 12 c_1 L_{\hat{g}}$,
  $\alpha \DefinedEqual \frac{4 c_1 L_{\hat{g}}}{\bar{M} - c_2 L_f}$,
  and any $\Gradient f(x^*) \in \Subdifferential f(x^*)$, it holds that
  \begin{multline*}
    \Expectation\Bigl[
      N [F(\tilde{x}_+) - F^*] + \frac{M_+}{2} \Norm{x_+ - x^*}^2
    \Bigr]
    \\
    \leq
    6 \alpha N
    \BregmanDistanceWithSubgradient{f}{\Gradient f(x^*)}(x^*, \tilde{x})
    +
    \frac{M}{2} \Norm{x - x^*}^2
    +
    N (c_3 \delta_f + 16 \alpha \delta_{\hat{g}})
    +
    c_4 D^2 \Expectation\bigl\{
      \PositivePart{\min\Set{M_+, \bar{M}} - M}
    \bigr\}.
  \end{multline*}
\end{lemma}

\begin{proof}
  Since $\hat{g}$ is an unbiased oracle for $\bar{g}$, so is $\hat{G}$
  (\cref{th:SvrgOracle}).
  Therefore, we can apply \cref{th:UniversalSgd-General} to get
  \begin{multline*}
    \Expectation\Bigl[
      N [F(\tilde{x}_+) - F^*]
      +
      \frac{M_+}{2} \Norm{x_+ - x^*}^2
      +
      \sum_{k = 0}^{N - 1} [
        \beta_{f, \bar{f}, \bar{g}}(x_{k + 1}, x_k)
        +
        \beta_{f, \bar{f}, \bar{g}}(x_k, x^*)
      ]
    \Bigr]
    \\
    \leq
    \frac{M}{2} \Norm{x - x^*}^2
    +
    \frac{c_1}{\bar{M} - c_2 L_f}
    \sum_{k = 0}^{N - 1} \Expectation[
      \Variance_{\hat{G}}(x_{k + 1}) + \Variance_{\hat{G}}(x_k)
    ]
    +
    c_3 N \delta_f
    \\
    +
    c_4 \Expectation\bigl\{
      \PositivePart{\min\Set{M_+, \bar{M}} - M} D^2
    \bigr\},
  \end{multline*}
  where $\bar{M} > c_2 L_f$ is an arbitrary constant and
  $x_k$ are the points generated inside $\UniversalSgd$.

  Applying now \cref{%
    th:SvrgOracle,%
    th:YoungsInequalityForVariance,%
    as:ApproximatelySmoothVariance,%
    as:ApproximatelySmoothVariance-Extra%
  }, we can estimate, for each $k$,
  \begin{align*}
    \hspace{2em}&\hspace{-2em}
    \Variance_{\hat{G}}(x_{k + 1})
    +
    \Variance_{\hat{G}}(x_k)
    =
    \Variance_{\hat{g}}(x_{k + 1}, \tilde{x})
    +
    \Variance_{\hat{g}}(x_k, \tilde{x})
    \leq
    2 \Variance_{\hat{g}}(x_{k + 1}, x_k)
    +
    3 \Variance_{\hat{g}}(x_k, \tilde{x})
    \\
    &\leq
    2 \Variance_{\hat{g}}(x_{k + 1}, x_k)
    +
    6 \Variance_{\hat{g}}(x_k, x^*)
    +
    6 \Variance_{\hat{g}}(x^*, \tilde{x})
    \\
    &\leq
    4 L_{\hat{g}} [\beta_{f, \bar{f}, \bar{g}}(x_{k + 1}, x_k) + \delta_{\hat{g}}]
    +
    12 L_{\hat{g}} [\beta_{f, \bar{f}, \bar{g}}(x_k, x^*) + \delta_{\hat{g}}]
    +
    24 L_{\hat{g}} [
      \BregmanDistanceWithSubgradient{f}{\Gradient f(x^*)}(x^*, \tilde{x})
      +
      2 \delta_{\hat{g}}
    ]
    \\
    &=
    4 L_{\hat{g}} [
      \beta_{f, \bar{f}, \bar{g}}(x_{k + 1}, x_k)
      +
      3 \beta_{f, \bar{f}, \bar{g}}(x_k, x^*)
      +
      6 \BregmanDistanceWithSubgradient{f}{\Gradient f(x^*)}(x^*, \tilde{x})
      +
      16 \delta_{\hat{g}}
    ],
  \end{align*}
  where $\Gradient f(x^*) \in \Subdifferential f(x^*)$ is arbitrary.
  Denoting
  $\alpha \DefinedEqual \frac{4 c_1 L_{\hat{g}}}{\bar{M} - c_2 L_f}$,
  we thus obtain
  \begin{multline*}
    \Expectation\Bigl[
      N [F(\tilde{x}_+) - F^*]
      +
      \frac{M_+}{2} \Norm{x_+ - x^*}^2
      +
      \sum_{k = 0}^{N - 1} [
        \beta_{f, \bar{f}, \bar{g}}(x_{k + 1}, x_k)
        +
        \beta_{f, \bar{f}, \bar{g}}(x_k, x^*)
      ]
    \Bigr]
    \\
    \leq
    6 \alpha N
    \BregmanDistanceWithSubgradient{f}{\Gradient f(x^*)}(x^*, \tilde{x})
    +
    \frac{M}{2} \Norm{x - x^*}^2
    +
    N (c_3 \delta_f + 16 \alpha \delta_{\hat{g}})
    +
    c_4 \Expectation\bigl\{
      \PositivePart{\min\Set{M_+, \bar{M}} - M} D^2
    \bigr\}
    \\
    +
    \alpha \sum_{k = 0}^{N - 1} \Expectation[
      \beta_{f, \bar{f}, \bar{g}}(x_{k + 1}, x_k)
      +
      3 \beta_{f, \bar{f}, \bar{g}}(x_k, x^*)
    ].
  \end{multline*}
  Requiring now $\bar{M} \geq c_2 L_f + 12 c_1 L_{\hat{g}}$,
  we get $\alpha \leq \frac{1}{3}$ which allows us to cancel the nonnegative
  $\beta_{f, \bar{f}, \bar{g}}(\cdot, \cdot)$ terms on both sides.
  The claim now follows.
\end{proof}

\thUniversalSvrg*

\begin{proof}
  \label{th:UniversalSvrg:Proof}
  The algorithm iterates
  $
    (\tilde{x}_{t + 1}, x_{t + 1}, M_{t + 1})
    \EqualRandom
    \UniversalSgd_{\hat{G}_t, \psi}(x_t, M_t, 2^{t + 1}; D)
  $
  for $t \geq 0$, where $\hat{G}_t = \SvrgOracle_{\hat{g}}(\tilde{x}_t)$.
  Applying \cref{th:UniversalSvrgEpoch}
  with $\bar{M} \DefinedEqual c_2 L_f + 48 c_1 L_{\hat{g}}$
  (for which $\alpha = \frac{1}{12}$ so that $6 \alpha 2^{t + 1} = 2^t$)
  and passing to full expectations, we obtain, for any $t \geq 0$,
  \begin{multline*}
    \Expectation\Bigl[
      2^{t + 1} [F(\tilde{x}_{t + 1}) - F^*]
      +
      \frac{M_{t + 1}}{2} \Norm{x_{t + 1} - x^*}^2
    \Bigr]
    \\
    \leq
    \Expectation\Bigl[
      2^t \beta_t
      +
      \frac{M_t}{2} \Norm{x_t - x^*}^2
      +
      c_4 \PositivePart{\min\Set{M_{t + 1}, \bar{M}} - M_t} D^2
    \Bigr]
    +
    2^{t + 1} \Bigl( c_3 \delta_f + \frac{4}{3} \delta_{\hat{g}} \Bigr),
  \end{multline*}
  where
  $
    \beta_t
    \DefinedEqual
    \BregmanDistanceWithSubgradient{f}{\Gradient f(x^*)}(x^*, \tilde{x}_t)
  $
  and $\Gradient f(x^*) \in \Subdifferential f(x^*)$ can be chosen arbitrarily.
  Rewriting $F_{t + 1} \DefinedEqual F(\tilde{x}_{t + 1}) - F^*$ as
  $F_{t + 1} = \beta_{t + 1} + (F_{t + 1} - \beta_{t + 1})$
  and telescoping the above inequalities
  (using, \cref{th:TelescopingDifferencesWithMin}),
  we get, for any $t \geq 1$,
  \begin{align*}
    \hspace{2em}&\hspace{-2em}
    \Expectation\Bigl[
      2^t \beta_t
      +
      \sum_{i = 1}^t 2^i (F_i - \beta_i)
      +
      \frac{M_t}{2} \Norm{x_t - x^*}^2
    \Bigr]
    \\
    &\leq
    \beta_0
    +
    \frac{M_0}{2} \Norm{x_0 - x^*}^2
    +
    \Bigl( c_3 \delta_f + \frac{4}{3} \delta_{\hat{g}} \Bigr)
    \sum_{i = 1}^t 2^i
    +
    c_4 \Expectation\bigl\{
      \PositivePart{\min\Set{M_t, \bar{M}} - M_0} D^2
    \bigr\}
    \\
    &\leq
    \beta_0
    +
    2 (2^t - 1) \Bigl( c_3 \delta_f + \frac{4}{3} \delta_{\hat{g}} \Bigr)
    +
    c_4 \bar{M} D^2
    \EqualDefines
    \Phi_0,
  \end{align*}
  where the final inequality is due to the fact that $M_0 = 0$, while
  $\sum_{i = 1}^t 2^i = 2 (2^t - 1)$.
  According to \cref{th:PropertyOfOptimalityCompositeOptimization},
  we can choose $\Gradient f(x^*) \in \Subdifferential f(x^*)$ such that
  $\beta_i \leq F(\tilde{x}_i) - F^*$ for all $i \geq 0$.
  Dropping now various nonnegative terms from the left-hand side of the above
  display, we conclude that
  \[
    2^t \Expectation[F_t] \leq \Phi_0.
  \]

  Let us estimate $\Phi_0$.
  Using our \cref{as:ApproximateSmoothness,as:BoundedFeasibleSet}
  and \cref{th:InequalitiesForApproximatelySmoothFunctions}
  (inequality~\eqref{eq:ApproximateSmoothness-UpperBoundOnUsualBregmanDistance}),
  we can bound
  $
    \beta_0
    \leq
    L_f \Norm{\tilde{x}_0 - x^*}^2 + 2 \delta_f
    \leq
    L_f D^2 + 2 \delta_f
  $.
  Therefore,
  \[
    \Phi_0
    \leq
    L_f D^2 + 2 \delta_f
    +
    c_4 \bar{M} D^2
    +
    2 (2^t - 1) (c_3 \delta_f + \tfrac{4}{3} \delta_{\hat{g}})
    \leq
    L D^2
    +
    2 (\bar{c}_3 \delta_f + \tfrac{4}{3} \delta_{\hat{g}}) \cdot 2^t
  \]
  where
  $
    L
    \DefinedEqual
    L_f + c_4 \bar{M}
    \equiv
    (c_2 c_4 + 1) L_f + 48 c_1 c_4 L_{\hat{g}}
  $
  and $\bar{c}_3 \DefinedEqual \max\Set{c_3, 1}$.
  Thus,
  \[
    \Expectation[F_t]
    \leq
    \frac{\Phi_0}{2^t}
    \leq
    \frac{L D^2}{2^t}
    +
    2 \bar{c}_3 \delta_f
    +
    \frac{8}{3} \delta_{\hat{g}},
  \]
  which proves the claimed convergence rate.

  Let us now estimate the number of oracle queries.
  At each iteration $t$, the algorithm first queries~$\bar{g}$ to construct
  the SVRG oracle~$\hat{G}_t$ (by precomputing $\bar{g}(\tilde{x}_t)$).
  All other queries are then done only to~$\hat{G}_t$ or, equivalently,
  to~$\hat{g}$ inside $\UniversalSgd_{\hat{G}_t, \psi}$ which is run
  for $N_{t + 1} = 2^{t + 1}$ iterations and thus requiring $\BigO(N_{t + 1})$
  queries to~$\hat{g}$.
  Summing up, after $T$ iterations, we obtain the total number of
  $\sum_{t = 1}^T \BigO(N_t) = \sum_{t = 1}^T \BigO(2^t) = \BigO(2^T)$
  queries to~$\hat{g}$, and $T$ queries to~$\bar{g}$.
\end{proof}

\subsubsection*{Helper Lemmas}

\begin{lemma}
  \label{th:PropertyOfOptimalityCompositeOptimization}
  Let $\Map{F}{\RealField^d}{\RealFieldPlusInfty}$ be the function
  $F(x) \DefinedEqual f(x) + \psi(x)$, where
  $\Map{f}{\RealField^d}{\RealField}$ is a convex function,
  and $\Map{\psi}{\RealField^d}{\RealFieldPlusInfty}$ is a proper closed convex
  function.
  Let $x^*$ be a minimizer of $F$ and let $F^* \DefinedEqual F(x^*)$.
  Then, there exists $\Gradient f(x^*) \in \Subdifferential f(x^*)$ such that,
  for any $x \in \EffectiveDomain \psi$,
  \[
    F(x) - F^*
    \geq
    \BregmanDistanceWithSubgradient{f}{\Gradient f(x^*)}(x^*, x).
  \]
\end{lemma}

\begin{proof}
  Since $x^*$ is a minimizer of~$F$, we have
  $
    0
    \in
    \Subdifferential F(x^*)
    =
    \Subdifferential{f}(x^*) + \Subdifferential \psi(x^*)
  $.
  In other words, there exists $\Gradient f(x^*) \in \Subdifferential f(x^*)$
  such that
  $
    \Gradient \psi(x^*)
    \DefinedEqual
    -\Gradient f(x^*)
    \in
    \Subdifferential \psi(x^*)
  $.
  Consequently, for any $x \in \EffectiveDomain \psi$,
  \begin{align*}
    F(x) - F^*
    &=
    f(x) - f(x^*)
    +
    [\psi(x) - \psi(x^*)]
    \geq
    f(x) - f(x^*)
    +
    \InnerProduct{\Gradient \psi(x^*)}{x - x^*}
    \\
    &=
    f(x) - f(x^*) - \InnerProduct{\Gradient f(x^*)}{x - x^*}
    =
    \BregmanDistanceWithSubgradient{f}{\Gradient f(x^*)}(x^*, x).
    \qedhere
  \end{align*}
\end{proof}

\subsection{Universal Fast SVRG}

\begin{lemma}[Universal Triangle SVRG Step]
  \label{th:UniversalTriangleSvrgStep}
  Consider problem~\eqref{eq:Problem} under
  \cref{as:BoundedFeasibleSet,as:ApproximateSmoothness,as:ApproximatelySmoothVariance}.
  Let $\tilde{x}, v \in \EffectiveDomain \psi$ be points,
  $M \geq 0$ and $A, a > 0$ be coefficients,
  $\hat{G} \DefinedEqual \SvrgOracle_{\hat{g}}(\tilde{x})$.
  Further, let, for $A_+ \DefinedEqual A + a$,
  \begin{gather*}
    x \DefinedEqual \frac{A \tilde{x} + a v}{A_+},
    \quad
    \hat{G}_x \EqualRandom \hat{G}(x),
    \quad
    \hat{v}_+ = \ProximalMap_{\psi}(v, \hat{G}_x, M / a),
    \quad
    \hat{x}_+ = \frac{A \tilde{x} + a \hat{v}_+}{A_+},
    \\
    \hat{G}_{x_+} \EqualRandom \hat{G}(\hat{x}_+),
    \quad
    \hat{M}_+
    =
    \frac{a^2}{A_+}
    M_+\Bigl(
      \frac{A_+}{a^2} M, \frac{a^2}{A_+^2} D^2, x, \hat{x}_+,
      \hat{G}_x, \hat{G}_{x_+}
    \Bigr).
  \end{gather*}
  Then, for
  $
    \bar{M}
    \DefinedEqual
    c_2 L_f \frac{a^2}{A_+} + 6 c_1 L_{\hat{g}} \frac{a^2}{A}
  $,
  it holds that
  \begin{multline*}
    \Expectation\Bigl[
      A_+ [F(\hat{x}_+) - F^*] + \frac{\hat{M}_+}{2} \Norm{\hat{v}_+ - x^*}^2
    \Bigr]
    \\
    \leq
    A [F(\tilde{x}) - F^*]
    +
    \frac{M}{2} \Norm{v - x^*}^2
    +
    c_4 D^2 \Expectation\bigl\{
      \PositivePart{\min\Set{\hat{M}_+, \bar{M}} - M}
    \bigr\}
    +
    c_3 A_+ \delta_f
    +
    \frac{5}{3} A \delta_{\hat{g}}.
  \end{multline*}
\end{lemma}

\begin{proof}
  Since $\hat{g}$ is an unbiased oracle for $\bar{g}$, so is $\hat{G}$
  (\cref{th:SvrgOracle}).
  Therefore, we can apply \cref{th:UniversalTriangleStep} to obtain
  \begin{multline*}
    \Expectation\Bigl[
      A_+ [F(\hat{x}_+) - F^*]
      +
      \frac{\hat{M}_+}{2} \Norm{\hat{v}_+ - x^*}^2
      +
      A_+ \beta_{f, \bar{f}, \bar{g}}(\hat{x}_+, x)
    \Bigr]
    +
    A \beta_{f, \bar{f}, \bar{g}}(x, \tilde{x})
    \\
    \leq
    A [F(\tilde{x}) - F^*]
    +
    \frac{M}{2} \Norm{v - x^*}^2
    +
    \frac{c_1 a^2}{\bar{M} - c_2 L_f \frac{a^2}{A_+}}
    \Expectation[\Variance_{\hat{G}}(\hat{x}_+) + \Variance_{\hat{G}}(x)]
    \\
    +
    c_3 A_+ \delta_f
    +
    c_4 \Expectation\bigl\{
      \PositivePart{\min\Set{\hat{M}_+, \bar{M}} - M} D^2
    \bigr\},
  \end{multline*}
  where $\bar{M} > c_2 L_f \frac{a^2}{A_+}$ is an arbitrary coefficient.
  Using \cref{%
    th:SvrgOracle,%
    th:YoungsInequalityForVariance,%
    as:ApproximatelySmoothVariance%
  }, we can further bound
  \begin{align*}
    \Variance_{\hat{G}}(\hat{x}_+) + \Variance_{\hat{G}}(x)
    &=
    \Variance_{\hat{g}}(\hat{x}_+, \tilde{x})
    +
    \Variance_{\hat{g}}(x, \tilde{x})
    \leq
    2 \Variance_{\hat{g}}(\hat{x}_+, x)
    +
    3 \Variance_{\hat{g}}(x, \tilde{x})
    \\
    &\leq
    2 L_{\hat{g}} [
      2 \beta_{f, \bar{f}, \bar{g}}(\hat{x}_+, x)
      +
      3 \beta_{f, \bar{f}, \bar{g}}(x, \tilde{x})
      +
      5 \delta_{\hat{g}}
    ].
  \end{align*}
  Denoting
  $
    \alpha
    \DefinedEqual
    \frac{2 c_1 L_{\hat{g}} a^2}{\bar{M} - c_2 L_f \frac{a^2}{A_+}}
  $,
  we thus obtain
  \begin{multline*}
    \Expectation\Bigl[
      A_+ [F(\hat{x}_+) - F^*]
      +
      \frac{\hat{M}_+}{2} \Norm{\hat{v}_+ - x^*}^2
      +
      A_+ \beta_{f, \bar{f}, \bar{g}}(\hat{x}_+, x)
    \Bigr]
    +
    A \beta_{f, \bar{f}, \bar{g}}(x, \tilde{x})
    \\
    \leq
    A [F(\tilde{x}) - F^*]
    +
    \frac{M}{2} \Norm{v - x^*}^2
    +
    c_3 A_+ \delta_f
    +
    5 \alpha \delta_{\hat{g}}
    +
    c_4 \Expectation\bigl\{
      \PositivePart{\min\Set{\hat{M}_+, \bar{M}} - M} D^2
    \bigr\}
    \\
    +
    2 \alpha \Expectation[\beta_{f, \bar{f}, \bar{g}}(\hat{x}_+, x)]
    +
    3 \alpha \beta_{f, \bar{f}, \bar{g}}(x, \tilde{x}).
  \end{multline*}
  Choosing now
  $
    \bar{M}
    =
    c_2 L_f \frac{a^2}{A_+} + 6 c_1 L_{\hat{g}} \frac{a^2}{A}
  $,
  we get $\alpha = \frac{1}{3} A \ (\leq \frac{1}{3} A_+)$, which allows us to
  drop the nonnegative $\beta_{f, \bar{f}, \bar{g}}(\cdot, \cdot)$ terms from
  both sides.
  The claim now follows.
\end{proof}

\begin{lemma}[Universal Triangle SVRG Epoch]
  \label{th:UniversalTriangleSvrgEpoch}
  Consider problem~\eqref{eq:Problem} under
  \cref{as:BoundedFeasibleSet,as:ApproximateSmoothness,as:ApproximatelySmoothVariance}.
  Let $\tilde{x}, v \in \EffectiveDomain \psi$ be points,
  $M \geq 0$ and $A, a > 0$ be coefficients, $N \geq 1$ be an integer, and let
  \[
    (\tilde{x}_+, v_+, M_+)
    \EqualRandom
    \UniversalTriangleSvrgEpoch_{\hat{g}, \psi}(\tilde{x}, v, M, A, a, N; D),
  \]
  as defined by \cref{alg:UniversalTriangleSvrgEpoch}.
  Then, for $A_+ \DefinedEqual A + a$ and
  $
    \bar{M}
    \DefinedEqual
    c_2 L_f \frac{a^2}{A_+} + 6 c_1 L_{\hat{g}} \frac{a^2}{A}
  $,
  it holds that
  \begin{multline*}
    \Expectation\Bigl[
      A_+ N [F(\tilde{x}_+) - F^*]
      +
      \frac{M_+}{2} \Norm{v_+ - x^*}^2
    \Bigr]
    \\
    \leq
    A N [F(\tilde{x}) - F^*]
    +
    \frac{M}{2} \Norm{v - x^*}^2
    +
    c_4 D^2 \Expectation\bigl\{
      \PositivePart{\min\Set{M_+, \bar{M}} - M}
    \bigr\}
    +
    N \Bigl( c_3 A_+ \delta_f + \frac{5}{3} A \delta_{\hat{g}} \Bigr).
  \end{multline*}
\end{lemma}

\begin{proof}
  Each iteration~$k$ of the algorithm, when conditioned on~$v_k$,
  follows the construction from \cref{th:UniversalTriangleSvrgStep}
  (with $v = v_k$, $M = M_k$, $A = A_k$, $a = a_{k + 1}$,
    $A_+ = A_{k + 1}$, $x = x_k$, $\hat{G}_x = G_{x_k}$,
    $\hat{v}_+ = v_{k + 1}$, $\hat{x}_+ = x_{k + 1}$,
    $\hat{G}_{x_+} = G_{x_{k + 1}}$, $\hat{M}_+ = M_{k + 1}$).
  Hence, we can write, after passing to full expectations,
  for each $k = 0, \ldots, N - 1$,
  \begin{multline*}
    \Expectation\Bigl[
      A_+ [F(x_{k + 1}) - F^*]
      +
      \frac{M_{k + 1}}{2} \Norm{v_{k + 1} - x^*}^2
    \Bigr]
    \\
    \leq
    A [F(\tilde{x}) - F^*]
    +
    \Expectation\Bigl[
      \frac{M_k}{2} \Norm{v_k - x^*}^2
      +
      c_4 \PositivePart{\min\Set{M_{k + 1}, \bar{M}} - M_k} D^2
    \Bigr]
    +
    \delta,
  \end{multline*}
  where
  $\delta \DefinedEqual c_3 A_+ \delta_f + \frac{5}{3} A \delta_{\hat{g}}$.
  Telescoping the above inequalities
  (using \cref{th:TelescopingDifferencesWithMin}), we get
  \begin{multline*}
    \Expectation\Bigl[
      A_+ \sum_{k = 1}^N [F(x_k) - F^*]
      +
      \frac{M_N}{2} \Norm{v_N - x^*}^2
    \Bigr]
    \\
    \leq
    A N [F(\tilde{x}) - F^*]
    +
    \frac{M_0}{2} \Norm{v_0 - x^*}^2
    +
    c_4 D^2 \Expectation\bigl\{
      \PositivePart{\min\Set{M_N, \bar{M}} - M_0}
    \bigr\}
    +
    N \delta.
  \end{multline*}
  The claim now follows from the convexity of~$F$ and our definitions
  $\tilde{x}_+ = \bar{x}_N = \frac{1}{N} \sum_{k = 1}^N x_k$,
  $v_+ = v_N$, $M_+ = M_N$, $M_0 = M$, $v_0 = v$.
\end{proof}

\ThUniversalFastSvrg*

\begin{proof}
  \label{th:UniversalFastSvrg:Proof}
  By our definition, the algorithm iterates for $t \geq 0$:
  \[
    (\tilde{x}_{t + 1}, v_{t + 1}, M_{t + 1})
    \EqualRandom
    \UniversalTriangleSvrgEpoch_{\hat{g}, \bar{g}, \psi}(
      \tilde{x}_t, v_t, M_t, A_t, a_{t + 1}, N; D
    ),
  \]
  where $A_t$ and $a_{t + 1}$ are deterministic coefficients satisfying the
  following equations:
  \begin{equation}
    \label{_2:EquationForCoefficients}
    A_{t + 1} = A_t + a_{t + 1}, \qquad a_{t + 1} = \sqrt{A_t}.
  \end{equation}
  In particular, for any $t \geq 0$, we have
  $
    \bar{M}_t'
    \DefinedEqual
    c_2 L_f \frac{a_{t + 1}^2}{A_{t + 1}}
    +
    6 c_1 L_{\hat{g}} \frac{a_{t + 1}^2}{A_t}
    \leq
    c_2 L_f + 6 c_1 L_{\hat{g}}
    \EqualDefines
    \bar{M}
  $,
  and hence
  $
    \PositivePart{\min\Set{M_{t + 1}, \bar{M}_t'} - M_t}
    \leq
    \PositivePart{\min\Set{M_{t + 1}, \bar{M}} - M_t}
  $
  (because, for any fixed~$a$ and~$b$, the function
    $\PositivePart{\min\Set{a, \cdot} - b}$ is nondecreasing as the composition
    of two nondecreasing functions).
  Applying now \cref{th:UniversalTriangleSvrgEpoch} and passing
  to full expectations, we therefore obtain, for any $t \geq 0$,
  \begin{multline*}
    \Expectation\Bigl[
      A_{t + 1} N [F(\tilde{x}_{t + 1}) - F^*]
      +
      \frac{M_{t + 1}}{2} \Norm{v_{t + 1} - x^*}^2
    \Bigr]
    \\
    \leq
    \Expectation\Bigl[
      A_t N [F(\tilde{x}_t) - F^*]
      +
      \frac{M_t}{2} \Norm{v_t - x^*}^2
      +
      c_4 \PositivePart{\min\Set{M_{t + 1}, \bar{M}} - M_t} D^2
    \Bigr]
    +
    N \Bigl( c_3 A_{t + 1} \delta_f + \frac{5}{3} A_t \delta_{\hat{g}} \Bigr).
  \end{multline*}

  Telescoping the above inequalities
  (using, in particular, \cref{th:TelescopingDifferencesWithMin}),
  we obtain, for any $t \geq 1$,
  \begin{multline*}
    A_t N \Expectation[F(\tilde{x}_t) - F^*]
    \leq
    A_0 N [F(\tilde{x}_0) - F^*]
    +
    \frac{M_0}{2} \Norm{v_0 - x^*}^2
    \\
    +
    c_4 \Expectation\bigl\{
      \PositivePart{\min\Set{M_t, \bar{M}} - M_0} D^2
    \bigr\}
    +
    N \Bigl(
      c_3 \delta_f \sum_{i = 1}^t A_i
      +
      \frac{5}{3} \delta_{\hat{g}} \sum_{i = 0}^{t - 1} A_i
    \Bigr)
    \\
    \leq
    A_0 N [F(\tilde{x}_0) - F^*]
    +
    c_4 \bar{M} D^2
    +
    N S_t (c_3 \delta_f + \tfrac{5}{3} \delta_{\hat{g}}),
  \end{multline*}
  where, for the last inequality, we have used the fact that
  $M_0 = 0$ and denoted $S_t \DefinedEqual \sum_{i = 1}^t A_i$.
  Thus, for any $t \geq 1$,
  \[
    \Expectation[F(\tilde{x}_t)] - F^*
    \leq
    \frac{1}{A_t} \Bigl(
      A_0 [F(\tilde{x}_0) - F^*] + \frac{c_4 \bar{M} D^2}{N}
    \Bigr)
    +
    \frac{S_t}{A_t} \Bigl( c_3 \delta_f + \frac{5}{3} \delta_{\hat{g}} \Bigr).
  \]

  At the same time, according to \cref{_2:EquationForCoefficients},
  $A_{t + 1} - A_t = \sqrt{A_t}$ for any $t \geq 0$.
  Hence, by \cref{th:SquaredGrowthWithDelay} (and our assumption on~$A_0$),
  we can estimate $A_t \geq \frac{1}{9} (t - t_0 + 1)^2$
  for any
  $
    t
    \geq
    t_0
    \DefinedEqual
    \Ceil{\log_2 \log_3 \frac{1}{A_0}} - 1
    \
    (\geq 0)
  $.
  Further, since the sequence~$A_t$ is increasing, we can
  estimate $S_t \equiv \sum_{i = 1}^t A_i \leq t A_t$,
  so that $\frac{S_t}{A_t} \leq t$.

  Substituting these bounds into the above display and using our formula
  for~$A_0$, we obtain, for any $t \geq t_0$,
  \[
    \Expectation[F(\tilde{x}_t)] - F^*
    \leq
    \rho_t [F(\tilde{x}_0) - F^* + c_4 \bar{M} D^2]
    +
    t (c_3 \delta_f + \tfrac{5}{3} \delta_{\hat{g}}),
  \]
  where $\rho_t \DefinedEqual \frac{9}{N (t - t_0 + 1)^2} \leq 1$.
  By our choice of~$\tilde{x}_0$, it holds that
  $F(\tilde{x}_0) - F^* \leq \frac{1}{2} L_f D^2 + \delta_f$
  (see \cref{th:FullGradientStepGivesGoodInitialPoint}).
  Denoting
  $
    L
    \DefinedEqual
    \frac{1}{2} L_f + c_4 \bar{M}
    \equiv
    (c_2 c_4 + \frac{1}{2}) L_f + 6 c_1 c_4 L_{\hat{g}}
  $,
  we get
  \[
    \Expectation[F(\tilde{x}_t)] - F^*
    \leq
    \rho_t (L D^2 + \delta_f)
    +
    t (c_3 \delta_f + \tfrac{5}{3} \delta_{\hat{g}})
    \leq
    \rho_t L D^2
    +
    (c_3 t + 1) \delta_f + \tfrac{5}{3} t \delta_{\hat{g}},
  \]
  which is exactly the claimed convergence rate bound.

  Let us now estimate the number of oracle queries.
  At the beginning, the algorithm makes $1$ query to~$\bar{g}$ to
  compute~$\tilde{x}_0$.
  All other queries to the oracles are then done, at each iteration~$t$, only
  inside the call to $\UniversalTriangleSvrgEpoch$
  (\cref{alg:UniversalTriangleSvrgEpoch}).
  Each such a call needs only one query to~$\bar{g}$ to construct the SVRG
  oracle~$\hat{G}$ (by precomputing $\bar{g}(\tilde{x})$), and
  $\BigO(N)$ queries to~$\hat{g}$ (which implements each query to~$\hat{G}$).
  Summing up, we get, after $t$ iterations, the total number of
  $\BigO(N t)$ queries to~$\hat{g}$ and $\BigO(t)$ queries to~$\bar{g}$.
\end{proof}

\subsubsection*{Helper Lemmas}

\begin{lemma}[c.f.~Lemma 1.1~in~\cite{grapiglia2017regularized}]
  \label{th:SquaredGrowthWithDelay}
  Let $A_t$ be a positive sequence such that
  \[
    A_{t + 1} - A_t \geq \sqrt{\gamma A_t}
  \]
  for all $t \geq 0$, where $\gamma > 0$, and let $A_0 \leq \frac{1}{9} \gamma$.
  Then, for any $t \geq 0$, we have
  \[
    A_t
    \geq
    \begin{cases}
      \gamma (\frac{A_0}{\gamma})^{1 / 2^t},
      &
      \text{if $t < t_0$},
      \\
      \frac{\gamma}{9} (t - t_0 + 1)^2,
      &
      \text{if $t \geq t_0$},
    \end{cases}
  \]
  where
  $
    t_0
    \DefinedEqual
    \Ceil{\log_2 \log_3 \frac{\gamma}{A_0}} - 1
    \ (\geq 0)
  $.
\end{lemma}

\begin{proof}
  By replacing $A_t$ with $A_t' = A_t / \gamma$, we can assume w.l.o.g.\ that
  $\gamma = 1$.

  For any $t \geq 0$, we have $A_{t + 1} \geq \sqrt{A_t}$, and hence
  \[
    A_t \geq A_0^{1 / 2^t}.
  \]
  In particular, for $t_0$ (as defined in the statement), we get
  $t_0 \geq \log_2 \log_3 \frac{1}{A_0} - 1$, so
  $2^{t_0} \geq \frac{1}{2} \log_3 \frac{1}{A_0}$, and hence
  \[
    A_{t_0}
    \geq
    A_0^{2 / \log_3 (1 / A_0)}
    =
    \bigl( 3^{-\log_3 (1 / A_0)} \bigr)^{2 / \log_3 (1 / A_0)}
    =
    3^{-2}
    =
    \frac{1}{9}
  \]
  (recall that $A_0 \leq \frac{1}{9} \leq 1$).

  On the other hand, for any $t \geq t_0$, we have
  \begin{align*}
    \sqrt{A_{t + 1}} - \sqrt{A_t}
    &\geq
    \sqrt{A_t + \sqrt{A_t}} - \sqrt{A_t}
    =
    \frac{\sqrt{A_t}}{\sqrt{A_t + \sqrt{A_t}} + \sqrt{A_t}}
    \\
    &=
    \frac{1}{\sqrt{1 + \frac{1}{\sqrt{A_t}}} + 1}
    \geq
    \frac{1}{\sqrt{1 + 3} + 1}
    =
    \frac{1}{3},
  \end{align*}
  where we have used the fact that $A_t \geq A_{t_0} \geq \frac{1}{9}$
  since $A_t$ is monotonically increasing.
  Telescoping these inequalities and rearranging, we get, for any $t \geq t_0$,
  \[
    A_t
    \geq
    \biggl( \frac{1}{3} (t - t_0) + \sqrt{A_{t_0}} \biggr)^2
    \geq
    \biggl( \frac{1}{3} (t - t_0) + \frac{1}{3} \biggr)^2
    =
    \frac{1}{9} (t - t_0 + 1)^2.
    \qedhere
  \]
\end{proof}

\begin{lemma}
  \label{th:FullGradientStepGivesGoodInitialPoint}
  Consider problem~\eqref{eq:Problem} under
  \cref{as:BoundedFeasibleSet,as:ApproximateSmoothness}.
  Let $x \in \EffectiveDomain \psi$, and let
  $x_+ \DefinedEqual \ProximalMap_{\psi}(x, \bar{g}(x), 0)$.
  Then, $F(x_+) - F^* \leq \frac{1}{2} L_f D^2 + \delta_f$.
\end{lemma}

\begin{proof}
  From the first-order optimality condition for the point~$x_+$
  (see \cref{th:OptimalityConditionForProximalStep}), it follows that
  \[
    \InnerProduct{\bar{g}(x)}{x^* - x_+} + \psi(x^*) \geq \psi(x_+).
  \]
  Combining the above inequality first with
  $
    f(x_+)
    \leq
    \bar{f}(x)
    +
    \InnerProduct{\bar{g}(x)}{x_+ - x}
    +
    \frac{L_f}{2} \Norm{x_+ - x}^2
    +
    \delta_f
  $
  and then with $\bar{f}(x) + \InnerProduct{\bar{g}(x)}{x^* - x} \leq f(x^*)$
  (which are both due to our \cref{as:ApproximateSmoothness}), we obtain
  \begin{align*}
    F(x_+)
    &=
    f(x_+) + \psi(x_+)
    \leq
    f(x_+) + \InnerProduct{\bar{g}(x)}{x^* - x_+} + \psi(x^*)
    \\
    &\leq
    \bar{f}(x)
    +
    \InnerProduct{\bar{g}(x)}{x^* - x}
    +
    \psi(x^*)
    +
    \frac{L_f}{2} \Norm{x_+ - x}^2
    +
    \delta_f
    \\
    &\leq
    F^* + \frac{L_f}{2} \Norm{x_+ - x}^2 + \delta_f.
  \end{align*}
  It remains to bound $\Norm{x_+ - x} \leq D$.
\end{proof}

\section{Omitted Proofs for \cref{sec:ApplicationToHolderSmoothProblems}}
\label{sec:ApplicationToHolderSmoothProblems-Proofs}

We start with the observation that for our specific example all our main
assumptions are satisfied.

\begin{remark}
  \label{rem:ConstantsForProblemWithHolderSmoothComponents}
  Under the setting from \cref{ex:ProblemWithHolderSmoothComponents},
  \cref{as:ApproximateSmoothness,as:ApproximatelySmoothVariance,as:ApproximatelySmoothVariance-Extra}
  are satisfied with
  $\bar{f} = f$,
  $
    \bar{g}(x)
    =
    \Gradient f(x)
    \DefinedEqual
    \Expectation_{\xi}[\Gradient f_{\xi}(x)]
  $,
  any $\delta_f, \delta_{\hat{g}} > 0$ and
  \[
    L_f
    =
    \biggl[
      \frac{1 - \nu}{2 (1 + \nu) \delta_f}
    \biggr]^{\frac{1 - \nu}{1 + \nu}}
    [H_f(\nu)]^{\frac{2}{1 + \nu}},
    \qquad
    L_{\hat{g}}
    =
    \frac{1}{b}
    \biggl[
      \frac{1 - \nu}{2 (1 + \nu) \delta_{\hat{g}}}
    \biggr]^{\frac{1 - \nu}{1 + \nu}}
    [H_{\max}(\nu)]^{\frac{2}{1 + \nu}}.
  \]
  Further, the oracle~$\hat{g}_b$ satisfies \cref{as:UniformlyBoundedVariance}
  with
  $
    \sigma_b^2
    \DefinedEqual
    \sup_{x \in \EffectiveDomain \psi} \Variance_{\hat{g}_b}(x)
    =
    \frac{1}{b} \sigma^2
  $,
  and
  $
    \sigma_{*, b}^2
    \DefinedEqual
    \Variance_{\hat{g}_b}(x^*)
    =
    \frac{1}{b} \sigma_*^2
  $.
\end{remark}

\begin{proof}
  For $b = 1$, this follows from
  \cref{%
    th:HolderGradientImpliesApproximateSmoothness,%
    th:StandardStochasticGradientOracleIsApproximatelySmooth%
  }
  and our definitions of~$\sigma^2$ and $\sigma_*$.
  The general case $b \geq 1$ follows from the fact that the standard
  mini-batching of size~$b$ reduces each of the variances
  $\Variance_{\hat{g}_1}(\cdot)$ and $\Variance_{\hat{g}_1}(\cdot, \cdot)$
  in~$b$ times.
\end{proof}

The following auxiliary result will be useful throughout this section:

\begin{lemma}
  \label{th:ConjugateToInversePowerFunction}
  Let $a, b, p > 0$ be real.
  Then,
  \[
    \min_{t > 0} \Bigl\{ \frac{a}{t^p} + b t \Bigr\}
    =
    (p + 1) a^{\frac{1}{p + 1}} \biggl( \frac{b}{p} \biggr)^{\frac{p}{p + 1}}.
  \]
\end{lemma}

\begin{proof}
  The expression inside the $\min$ is a convex function in~$t > 0$.
  Differentiating and setting its derivative to zero, we see that the
  minimum is attained at $t_* = (\frac{a p}{b})^{\frac{1}{p + 1}}$.
  Hence,
  \[
    \min_{t > 0} \Bigl\{ \frac{a}{t^p} + b t \Bigr\}
    =
    a \biggl( \frac{b}{a p} \biggr)^{\frac{p}{p + 1}}
    +
    b \biggl( \frac{a p}{b} \biggr)^{\frac{1}{p + 1}}
    =
    (p + 1) a^{\frac{1}{p + 1}} \biggl( \frac{b}{p} \biggr)^{\frac{p}{p + 1}}.
    \qedhere
  \]
\end{proof}

\section{Additional Discussion of Related Work}

\paragraph{Inexact Oracle and Approximate Smoothness.}

\textcite{devolder2013first} introduced the notion of the inexact first-order
oracle and analyzed the behaviour of several first-order methods for smooth
convex optimization using such an oracle.
Although their work was motivated by the desire to present the general
definition of an inexact oracle covering many different applications,
it was also observed that this oracle model is suitable for studying
weakly smooth problems.
This insight was later used in~\cite{nesterov2015universal} to develop
universal gradient methods for Hölder smooth problems.
First stochastic gradient methods for approximately smooth functions with
inexact oracle were proposed in~\cite{devolder2011stochastic}.
These algorithms however are not adaptive and require the knowledge of
problem-dependent constants.
For more details on the subject, see~\cite{Devolder-13-PhdThesis}.

\paragraph{Parameter-Free Methods.}

Parameter-free algorithms originating from the literature on online
learning~\cite{%
  streeter2012noregret,%
  orabona2014simultaneous,%
  cutkosky2017online,%
  cutkosky2018black,%
  mhammedi2020lipschitz,%
  jacobsen2023unconstrained%
}
is another popular type of adaptive methods.
They are usually endowed with mechanisms helping achieving efficiency bounds
that are almost insensitive (typically, with logarithmic dependency) to the
error of estimating certain problem parameters, such as the diameter of the
feasible set~\cite{%
  carmon2022making,%
  defazio2023learning,%
  ivgi2023dog,%
  khaled2023dowg,%
  mishchenko2023prodigy%
}.

\paragraph{Variance Reduction.}

Variance reduction techniques encompass a set of strategies that enhance the
convergence speed of SGD when multiple passes are possible over the training
dataset.
Various researchers simultaneously introduced methods to reduce variance around
the same period~\cite{%
  johnson2013accelerating,%
  zhang2013linear,%
  mahdavi2013mixed,%
  wang2013variance%
}.
The consideration of mini-batching in the context of these methods
is documented in~\cite{babanezhad2015stopwasting}, while,
in~\cite{gower2021stochastic}, it is shown that the convergence rate is
influenced by both the average and the maximum smoothness of individual
components.
For further details, see~\cite{gower2020variance} and the references therein.

Sometimes, it is even not necessary to use an explicit variance reduction
mechanism.
SGD may converge fast in the so-called over-parameterized regime, or when the
stochastic noise is small at the optimal solution~\cite{%
  cotter2011better,%
  schmidt2013fast,%
  needell2014stochastic,%
  ma2018power,%
  liu2018accelerating,%
  necoara2019linear%
}.
In this work, we call this effect implicit variance reduction.
Such a situation is also considered
in~\cite{stich2019unified,gorbunov2020unified} and, more recently,
\textcite{woodworth2021even} proposed an accelerated SGD algorithm for this
setting, under the assumption that the smoothness and noise constants are
known.

\section{Additional Experiments}
\label{sec:AdditionalExperiments}

\subsection{Logistic Regression with Real-World Data}
\label{sec:AdditionalExperiments:LogisticRegression}

In this section, we present experiments on the
\emph{logistic regression problem}:
\[
  f^*
  =
  \min_{\Norm{x} \leq R} \Bigl\{
    f(x)
    \DefinedEqual
    \frac{1}{n} \sum_{i = 1}^n \log(1 + e^{-b_i \InnerProduct{a_i}{x}})
  \Bigr\},
\]
where $a_i \in \RealField^d$ and $b_i \in \Set{-1, 1}$ are features and labels
taken from diverse real-world datasets from LIBSVM~\cite{chang2011libsvm}:
mushrooms ($d \ll n$), w8a ($d \ll n$), leu ($d \gg n$) and
colon-cancer ($d \gg n$).
The dataset leu is quite special because it satisfies the so-called
interpolation condition, meaning that the variance at the optimum is zero.
We fix $R = 1$ and use the mini-batch size of~$b = 32$ for the first two
datasets and $b = 1$ for the last two.

\begin{figure}[tb]
  \centering
  \includegraphics[width=\textwidth]{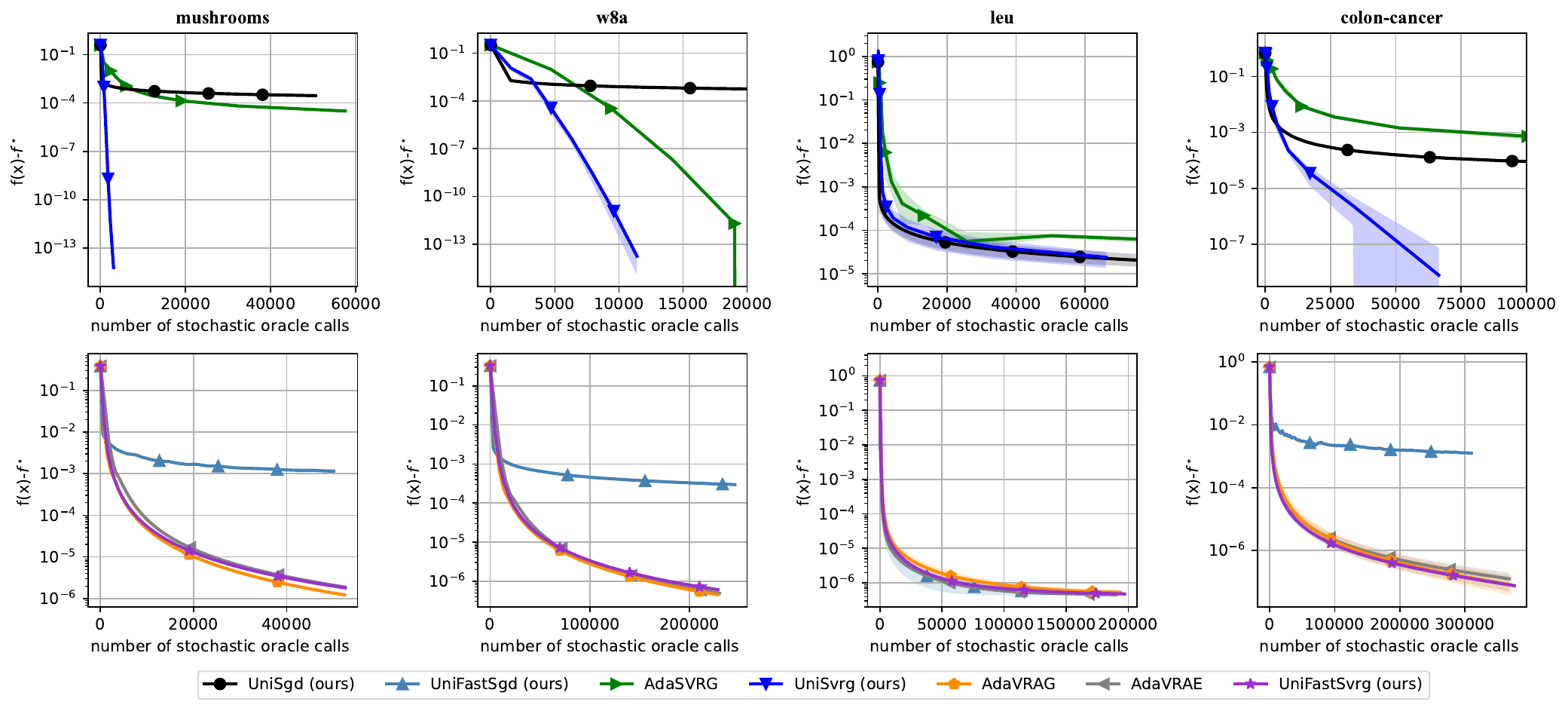}
  \caption{
    Comparison of various methods on the logistic regression problem with
    real-world data.
  }
  \label{fig:LogisticRegression}
\end{figure}

\Cref{fig:LogisticRegression} shows the results of our experiments.
The solid lines and the shaded area for each method represent, respectively,
the mean and the region between the minimum and the maximum values after three
independent runs of the algorithm.
We see that, on the leu dataset, $\UniversalSgd$ and $\UniversalFastSgd$
converge as fast as the best non-accelerated and accelerated SVRG methods,
respectively, which confirms our theory on implicit variance reduction.
Otherwise, these two SGD methods are typically much slower than the SVRG
algorithms.
Our $\UniversalSvrg$ method performs consistently better than AdaSVRG across
all the datasets.
Overall, all adaptive accelerated SVRG methods demonstrate comparable
performance for solving these smooth problems.

\subsection{Comparison between Stepsize Update Rules}
\label{sec:AdditionalExperiments:ComparisonBetweenStepsizeUpdateRules}

In this section, we compare the AdaGrad stepsize
rule~\eqref{eq:AdaGradStepsizeUpdateRule} with the other
rule~\eqref{eq:ModifiedAdaGradRule}
for $\UniversalSgd$ (\cref{alg:UniversalSgd}),
$\UniversalFastSgd$ (\cref{alg:UniversalFastSgd}),
$\UniversalSvrg$ (\cref{alg:UniversalSvrg}), and
$\UniversalFastSvrg$ (\cref{alg:UniversalFastSvrg}).
We consider the polyhedron feasibility and logistic regression problems
under the same setups as in
\cref{sec:Experiments,sec:AdditionalExperiments:LogisticRegression}.

The results are shown in \cref{%
  fig:CompareStepsizeUpdateRules-PolyhedronFeasibility,%
  fig:CompareStepsizeUpdateRules-LogisticRegression%
},
where we plot the function residual and the stepsize (inverse of $M$)
against stochastic oracle calls.
We see that the two stepsize rules work very similarly across all test cases,
which was not evident from the theory alone.

\begin{figure}[tb!]
  \centering
  \includegraphics[width=\textwidth]{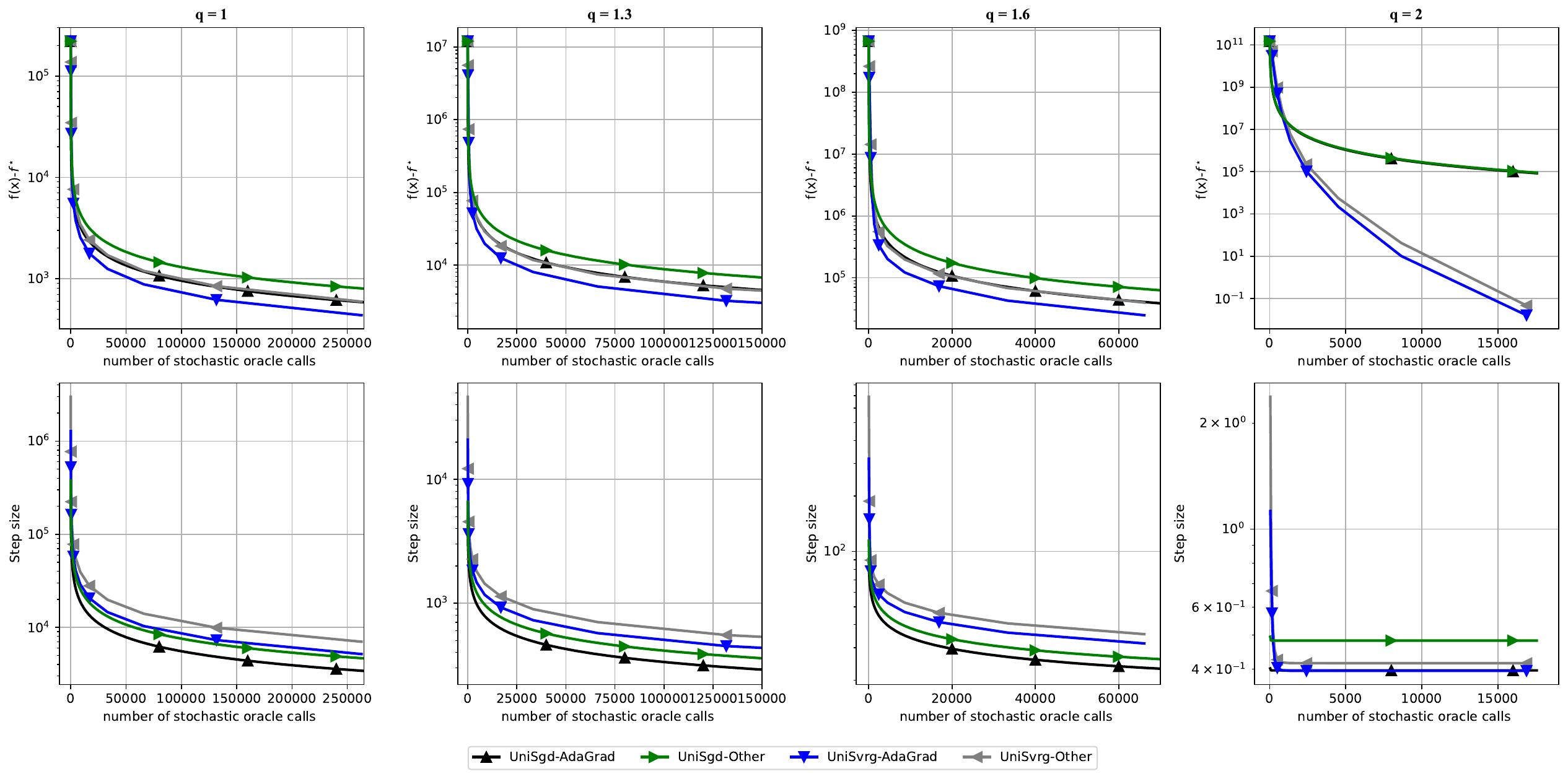}
  \includegraphics[width=\textwidth]{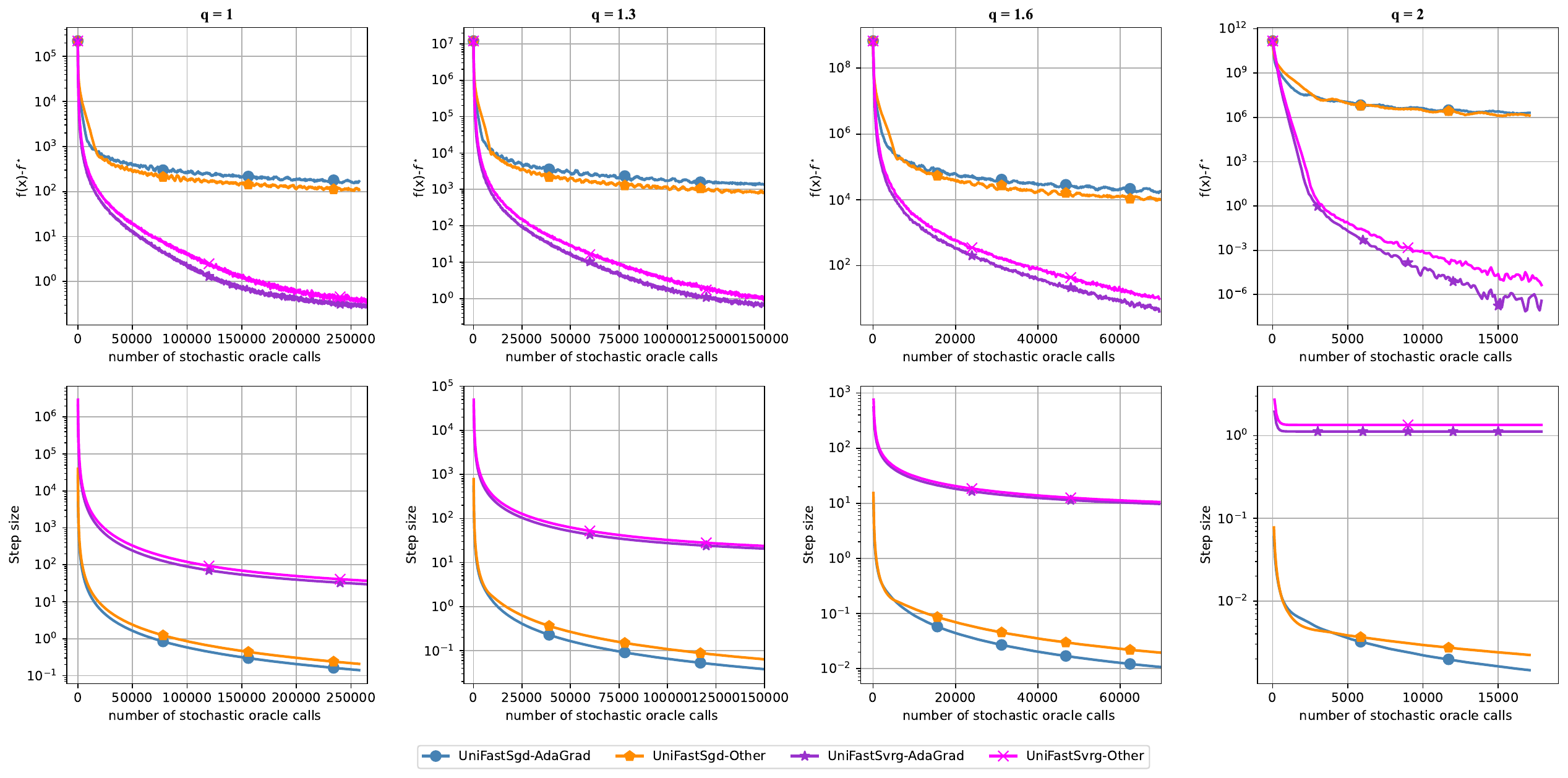}
  \caption{
    Comparison of our methods for different stepsize update rules on the
    polyhedron feasibility problem.
  }
  \label{fig:CompareStepsizeUpdateRules-PolyhedronFeasibility}
\end{figure}

\begin{figure}[tb]
  \centering
  \includegraphics[width=\textwidth]{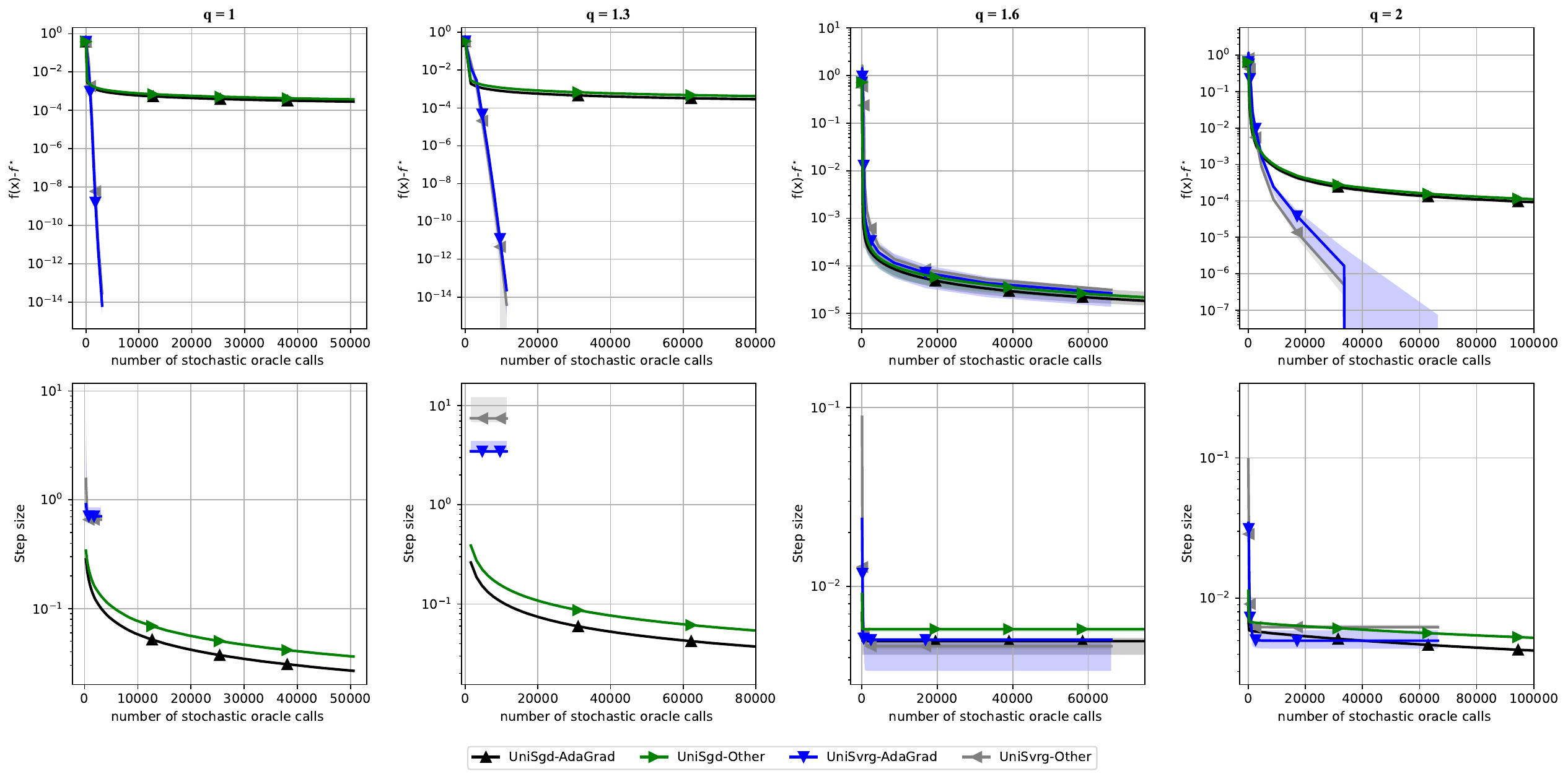}
  \includegraphics[width=\textwidth]{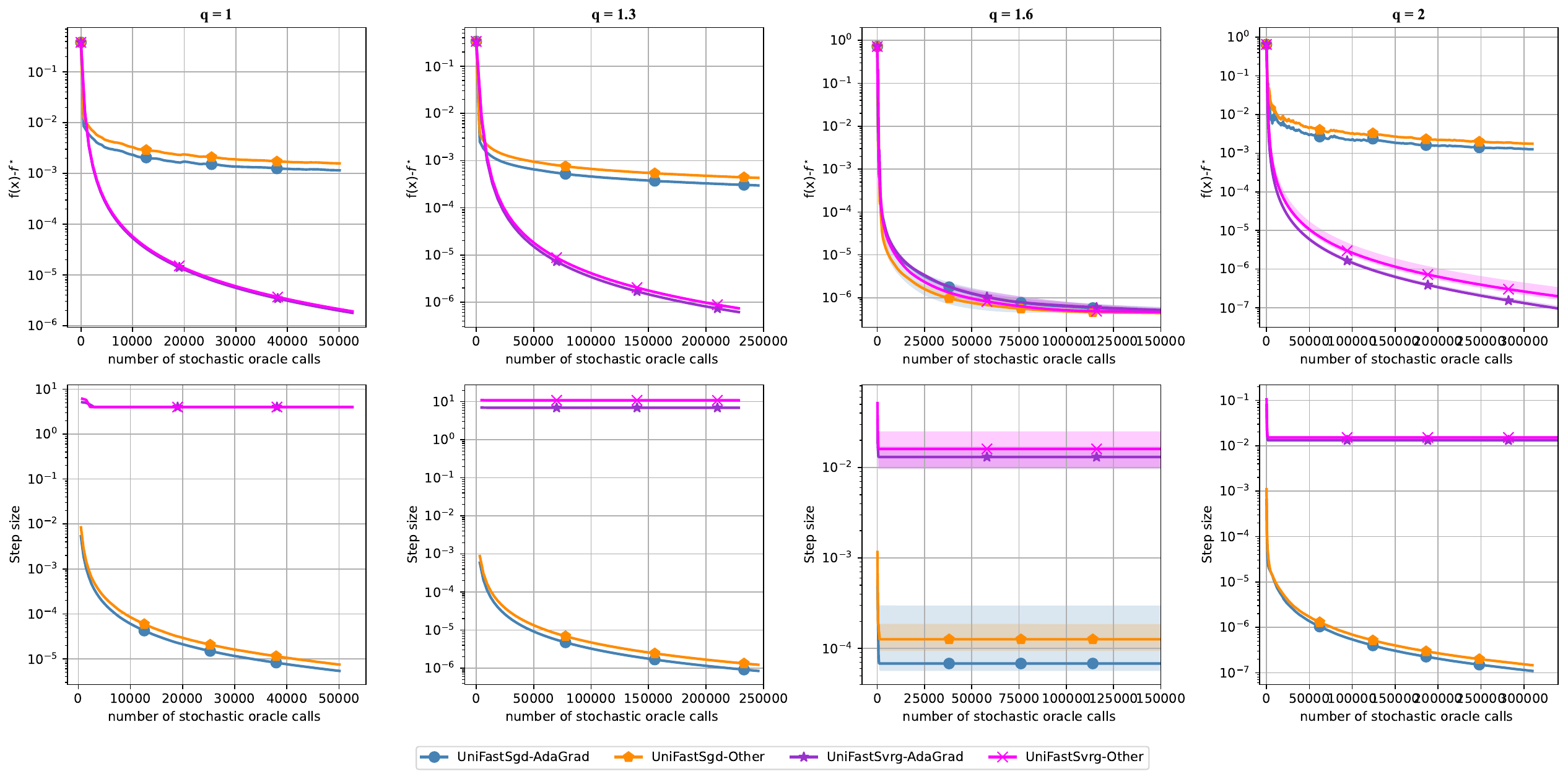}
  \caption{
    Comparison of our methods for different stepsize update rules on the
    logistic regression problem with real-world data.
  }
  \label{fig:CompareStepsizeUpdateRules-LogisticRegression}
\end{figure}

\end{document}